\documentclass[12pt]{amsart}
\usepackage{graphicx, overpic, epsf,amssymb,amscd,amsfonts}

\newtheorem{thm}{Theorem}[section]
\newtheorem*{claim}{Claim}

\newtheorem{prop}[thm]{Proposition}

\newtheorem{lemma}[thm]{Lemma}
\newtheorem{cor}[thm]{Corollary}

\newcommand{\C}{\mathbb{C}}
\newcommand{\R}{\mathbb{R}}
\newcommand{\Z}{\mathbb{Z}}

\newcommand{\bdry}{\partial}

\newcommand{\s}{\vskip.1in}
\newcommand{\n}{\noindent}

\newcommand{\be}{\begin{enumerate}}
\newcommand{\ee}{\end{enumerate}}

\numberwithin{equation}{subsection}

\def\SFH{\textit{SFH} \hspace{.03in}}
\def\CF{\textit{CF} \hspace{.03in}}
\def\EH{\textit{EH} \hspace{.03in}}
\def\HFhat{\widehat{\textit{HF}} \hspace{.03in}}
\def\HFKhat{\widehat{\textit{HFK}} \hspace{.03in}}

\begin{document}
%%%%%%%%%%%%%%%%%%%%%%%%%%%%%%%%%%%%%%%%

\title{The contact invariant in sutured Floer homology}

\author{Ko Honda}
\address{University of Southern California, Los Angeles, CA 90089}
\email{khonda@usc.edu} \urladdr{http://rcf.usc.edu/\char126 khonda}

\author{William H. Kazez}
\address{University of Georgia, Athens, GA 30602} \email{will@math.uga.edu}
\urladdr{http://www.math.uga.edu/\char126 will}

\author{Gordana Mati\'c}
\address{University of Georgia, Athens, GA 30602} \email{gordana@math.uga.edu}
\urladdr{http://www.math.uga.edu/\char126 gordana}

\date{This version: September 10, 2007. (The pictures are in color.)}



\thanks{KH supported by an NSF
CAREER Award (DMS-0237386); GM supported by NSF grant DMS-0410066;
WHK supported by NSF grant DMS-0406158.}

\begin{abstract}
We describe an invariant of a contact 3-manifold with convex
boundary as an element of Juh\'asz's sutured Floer homology. Our
invariant generalizes the contact invariant in Heegaard Floer
homology in the closed case, due to Ozsv\'ath and Szab\'o.
\end{abstract}

\maketitle

The goal of this paper is to present an invariant of a contact
3-manifold with convex boundary.  The contact class is an element of
{\em sutured Floer homology}, defined by Andr\'as Juh\'asz, and
generalizes the contact class in Heegaard Floer homology in the
closed case, as defined by Ozsv\'ath and Szab\'o~\cite{OS3} and
reformulated by the authors in~\cite{HKM2}. In this paper we assume
familiarity with convex surface theory (cf.~\cite{Gi1,H1}), Heegaard
Floer homology (cf.~\cite{OS1,OS2}), and sutured manifold theory
(cf.~\cite{Ga}).

Sutured manifold theory was introduced by Gabai~\cite{Ga} to study
and construct taut foliations. A {\em sutured manifold} $(M,\Gamma)$
is a compact oriented 3-manifold (not necessarily connected) with
boundary, together with a compact subsurface $\Gamma=A(\Gamma)\sqcup
T(\Gamma)\subset \bdry M$, where $A(\Gamma)$ is a union of pairwise
disjoint annuli and $T(\Gamma)$ is a union of tori.  We orient each
component of $\bdry M- \Gamma$, subject to the condition that the
orientation changes every time we nontrivially cross $A(\Gamma)$.
Let $R_+(\Gamma)$ (resp.\ $R_-(\Gamma)$) be the open subsurface of
$\bdry M-\Gamma$ on which the orientation agrees with (resp.\ is the
opposite of) the boundary orientation on $\bdry M$. Moreover,
$\Gamma$ is oriented so that the orientation agrees with the
orientation of $\bdry R_+(\Gamma)$ and is the opposite that of
$\bdry R_-(\Gamma)$. A sutured manifold $(M,\Gamma)$ is {\em
balanced} if $M$ has no closed components,
$\pi_0(A(\Gamma))\rightarrow \pi_0(\bdry M)$ is surjective, and
$\chi(R_+(\Gamma))=\chi(R_-(\Gamma))$ on every component of $M$. In
particular, $\Gamma=A(\Gamma)$ and every boundary component of
$\bdry M$ nontrivially intersects the suture $\Gamma$. {\em In this
paper, all our sutured manifolds are assumed to be balanced.}

In what follows, we view the suture $\Gamma$ as the union of cores
of the annuli, i.e., as a disjoint union of simple closed curves.

In the setting of contact structures, a natural condition to impose
on a contact 3-manifold $(M,\xi)$ with boundary is to require that
$\bdry M$ be {\em convex}, i.e., there is a contact vector field
transverse to $\bdry M$.  To a convex surface $F$ one can associate
its {\em dividing set} $\Gamma_F$ --- it is defined as the isotopy
class of multicurves $\{x\in F| X(x)\in \xi(x)\}$, where $X$ is a
transverse contact vector field. It was explained in \cite{HKM3}
that dividing sets and sutures on balanced sutured manifolds can be
viewed as equivalent objects. As usual, our contact structures are
assumed to be cooriented.

In a pair of important papers~\cite{Ju1,Ju2}, Andr\'as Juh\'asz
generalized (the hat versions of) Ozsv\'ath and Szab\'o's Heegaard
Floer homology~\cite{OS1,OS2} and link Floer homology~\cite{OS4}
theories, and assigned a Floer homology group $\SFH(M,\Gamma)$ to a
balanced sutured manifold $(M,\Gamma)$. (A similar theory was also
worked out by Lipshitz~\cite{Li2}.) An important property of this
{\em sutured Floer homology} is the following: if
$(M,\Gamma)\stackrel{T}\rightsquigarrow (M',\Gamma')$ is a sutured
manifold decomposition along a cutting surface $T$, then
$\SFH(M',\Gamma')$ is a direct summand of $\SFH(M,\Gamma)$.

The main result of this paper is the following:

\begin{thm}
Let $(M,\Gamma)$ be a balanced sutured manifold, and let $\xi$ be a
contact structure on $M$ with convex boundary, whose dividing set on
$\bdry M$ is $\Gamma$. Then there exists an invariant
$\EH(M,\Gamma,\xi)$ of the contact structure which lives in
$\SFH(-M,-\Gamma)/\{\pm 1\}$.
\end{thm}

Here we are using $\Z$-coefficients.  Note that there is currently a
$\pm 1$ ambiguity when $\Z$-coefficients are used.

The paper is organized as follows.  Section~\ref{section: prelim} is
devoted to discussing the contact-topological preliminaries for
obtaining a partial open book decomposition of a contact 3-manifold
with convex boundary.   We define the contact invariant in
Section~\ref{section: defn} and prove that it is independent of the
choices made in Section~\ref{section: well-defined}.  We discuss
some basic properties of the contact class in Section~\ref{section:
properties} and compute some examples in Section~\ref{section:
examples}.  Finally, we explain the relationship to sutured manifold
decompositions in Section~\ref{section: sutured}.

\section{Contact structure preliminaries} \label{section: prelim}

Let $(M,\Gamma)$ be a sutured manifold.  Let $\xi$ be a contact
structure on $M$ with convex boundary so that the dividing set
$\Gamma_{\bdry M}$ on $\bdry M$ is isotopic to $\Gamma$. Such a
contact manifold will be denoted $(M,\Gamma,\xi)$.

The following theorem is the key to obtaining a partial open book
decomposition, slightly generalizing the work of Giroux~\cite{Gi2}
to the relative case.  For more detailed expositions of Giroux's
work, see \cite{Co2,Et}.

\begin{thm}\label{thm: decomposition}
There exists a Legendrian graph $K\subset M$ whose endpoints, i.e.,
univalent vertices, lie on $\Gamma\subset \bdry M$ and which
satisfies the following: \be \item There is a neighborhood
$N(K)\subset M$ of $K$ so that (i) $\bdry N(K)= T\cup (\cup_i D_i)$,
(ii) $T$ is a convex surface with Legendrian boundary, (iii)
$D_i\subset \bdry M$ is  a convex disk with Legendrian boundary,
(iv) $T \cap \bdry M = \cup_i\bdry D_i$, (v) $\# (\bdry D_i\cap
\Gamma_{\bdry M})=2$, and (vi) there is a system of pairwise
disjoint compressing disks $D'_j$ for $N(K)$ so that $\bdry
D'_j\subset T$, $|\bdry D'_j\cap \Gamma_T|=2$, and each component of
$N(K)-\cup_j D'_j$ is a standard contact 3-ball, after rounding the
corners.
\item Each component $H$ of the complement $M-N(K)$ is a
handlebody with convex boundary. There is a system of pairwise
disjoint compressing disks $D^\alpha_k$ for $H$ so that $|\bdry
D^\alpha_k\cap \Gamma_{\bdry H}|=2$ and $H-\cup_k D^\alpha_k$ is a
standard contact 3-ball, after rounding the corners. \ee
\end{thm}

Here $|\cdot|$ denotes the geometric intersection number and
$\#(\cdot)$ denotes the number of connected components. A {\em
standard contact 3-ball} is a tight contact 3-ball $B^3$ with convex
boundary and $\#\Gamma_{\bdry B^3}=1$. We say that a handlebody $H$
with convex boundary admits a {\em product disk decomposition} if
condition (2) of Theorem~\ref{thm: decomposition} holds.

\begin{proof}
Since $F=\bdry M$ is convex, there is an $I=[0,1]$-invariant contact
neighborhood $F\times I\subset M$ so that $F\times\{1\}=\bdry M$.
First we take a polyhedral decomposition of $F_0=F\times\{0\}$ so
that the 1-skeleton is Legendrian and the boundary of each 2-cell
intersects $\Gamma_{F_0}$ in two points. During this process we need
to use the Legendrian realization principle and slightly isotop
$F_0$. Next, extend the polyhedral decomposition on $F_0$ to
$M-(F\times I)$ so that the 1-skeleton $K'$ is Legendrian and the
boundary of each 2-cell has Thurston-Bennequin invariant $tb=-1$.
Finally, we need to connect $K'$ to $\bdry M$. To achieve this, for
each intersection point $p$ of $K'|_{F_0}$ with $\Gamma_{F_0}$, we
add an edge $\{p\}\times[0,1]\subset F\times[0,1]$ to $K'$.  The
resulting graph is our desired $K$.

We now prove that $K$ satisfies the properties of the theorem. There
is a collection of compressing disks in $(M - N(K'))-(F\times I)$
which intersect the dividing set $\Gamma_{\bdry N(K')}$ at exactly
two points, by the $tb=-1$ condition.  Without loss of generality,
these compressing disks can be made convex with Legendrian boundary.
Then the compressing disks cut up $M-N(K')$ into a disjoint union of
standard contact 3-balls and $F\times I$. After removing standard
neighborhoods of $\{p\}\times I$ from $F\times I$, the remaining
contact 3-manifold admits a product disk decomposition.   This
implies that $M-N(K)$ admits a product disk decomposition.
\end{proof}

The next theorem is the relative version of the subdivision theorem
of Giroux~\cite{Gi2} for contact cellular decompositions, which
implies that on a closed manifold $M$ any two open books
corresponding to a fixed $(M,\xi)$ become isotopic after a sequence
of positive stabilizations to each.

\begin{thm}\label{thm: subdivision}
Let $K$ and $K'$ be Legendrian graphs on $(M,\Gamma,\xi)$ which
satisfy the conditions of Theorem~\ref{thm: decomposition} and
$K\cap K'=\emptyset$. Then there exists a common Legendrian
extension $L$ of $K$ and $K'$ so that $L=K_n$ is obtained
inductively from $K=K_0$ by attaching Legendrian arcs $c_i$,
$i=1,\dots,n-1$, to the standard neighborhood $N(K_i)$ of $K_i$ so
that the following hold:
\begin{enumerate}
\item The arc $c_i$ has both endpoints on the dividing set of $\bdry (M-N(K_i))$,
and $int(c_i)\subset int(M-N(K_i))$.
\item $N(K_{i+1})=N(c_i)\cup N(K_i)$, and $K_{i+1}$ is a Legendrian
graph so that $N(K_{i+1})$ is its standard neighborhood.
\item There is a Legendrian arc $d_i$ on $\bdry (M-N(K_i))$ with the same
endpoints as $c_i$, after possible application of the Legendrian
realization principle. The arc $d_i$ intersects $\Gamma_{\bdry
(M-N(K_i))}$ only at its endpoints.
\item \label{unknot}
The Legendrian knot $\gamma_i=c_i\cup d_i$ bounds a disk in
$M-N(K_i)$ and has $tb(\gamma_i)=-1$ with respect to this disk. This
implies that $c_i$ and $d_i$ are Legendrian isotopic relative to their
endpoints inside the closure of $M-N(K_i)$.
\end{enumerate}
The graph $L=K'_m$ is similarly obtained from $K'=K'_0$ by
inductively attaching Legendrian arcs $c_i'$ to $N(K_i')$ to obtain
$K_{i+1}'$.
\end{thm}

\begin{figure}[ht]
\begin{overpic}[width=8cm]{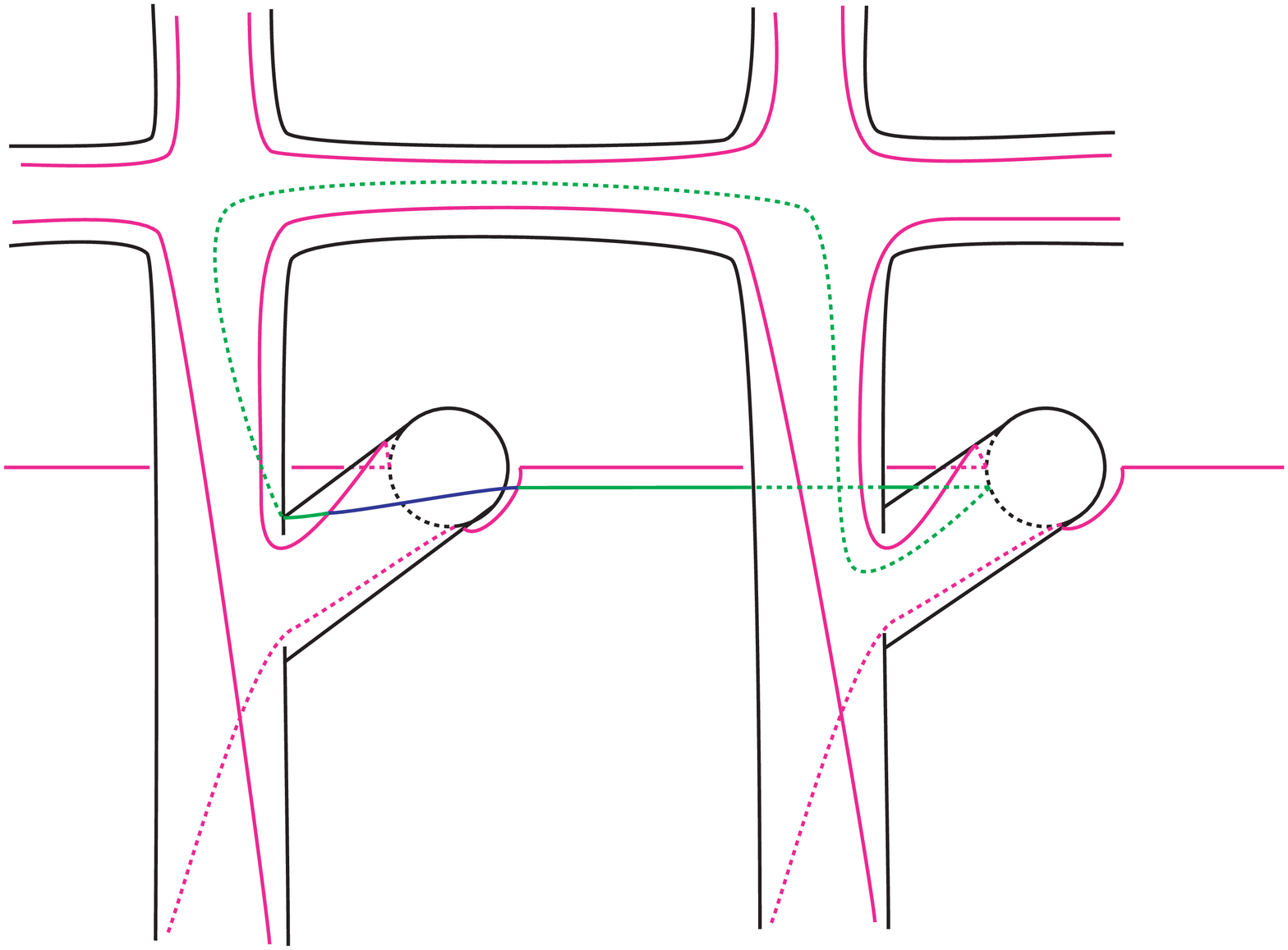}
\put(33,36){\tiny $c_i$} \put(47,32.3){\tiny $d_i$}
\end{overpic}
\caption{Step (3) of the subdivision process. The blue arc is $c_i$
and the green arc is $d_i$.} \label{attacharc2}
\end{figure}

\begin{proof}
Consider the invariant neighborhood $F\times[0,1]$ with $F=F_1=\bdry
M$ as in the first paragraph of Theorem~\ref{thm: decomposition},
and take a copy $F_{1-\delta}$ near $F$.  The procedure for finding
a common refinement $L$ of $K$ and $K'$ is as follows: \be
\item First add Legendrian arcs to $K$ so that $K_i$, $i\gg 0$,
contains a Legendrian $1$-skeleton of the protective layer
$F_{1-\delta}$ (as in the first paragraph of Theorem~\ref{thm:
decomposition}).
\item Subdivide the Legendrian $1$-skeleton of $F_{1-\delta}$
sufficiently, so that every arc of $K'$ near $\bdry M$ (we assume
they are all of the form $\{p\}\times [1-\delta,1]$) intersects the
$1$-skeleton.
\item Next add enough Legendrian arcs of the type $\{p\}\times
[1-\delta,1]$, where $p\in \Gamma_F$, so that $K_i\supset (K\cup
K')\cap (F\times[1-\delta,1])$, $i\gg 0$.
\item Finally, apply the contact subdivision procedure of \cite{Gi2}
away from $F\times[1-\delta,1]$. \ee

Steps (2) and (4) involve the same procedure used in \cite{Gi2},
namely, given a face $\Delta$ of the contact cellular decomposition
with Legendrian boundary $\bdry \Delta\subset K$, take a Legendrian
arc $c_i\subset \Delta$ with endpoints on $\bdry \Delta$, so that
$c_i$ cuts $\Delta$ into $\Delta_1$, $\Delta_2$, each of which has
$tb(\bdry \Delta_i)=-1$.

Next we discuss what happens in Steps (1) and (3).  There are two
types of arcs to attach in Step (1). The first type is an arc
$c_i\subset F_{1-\delta}$ from $(p,1-\delta)$ to $(q,1-\delta)$,
$p,q\in \Gamma_F$, which does not intersect
$\Gamma_F\times\{1-\delta\}$ in its interior.
Figure~\ref{attacharc1} depicts this situation.
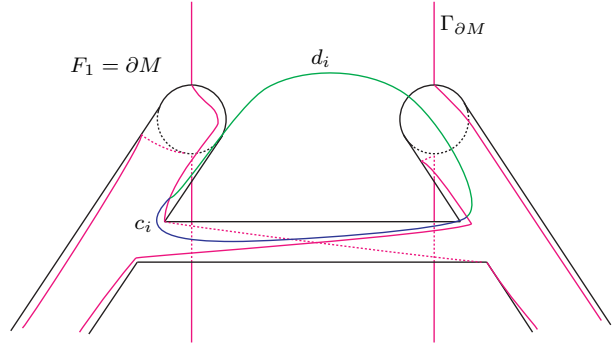
\begin{figure}[ht]
\begin{overpic}[width=8cm]{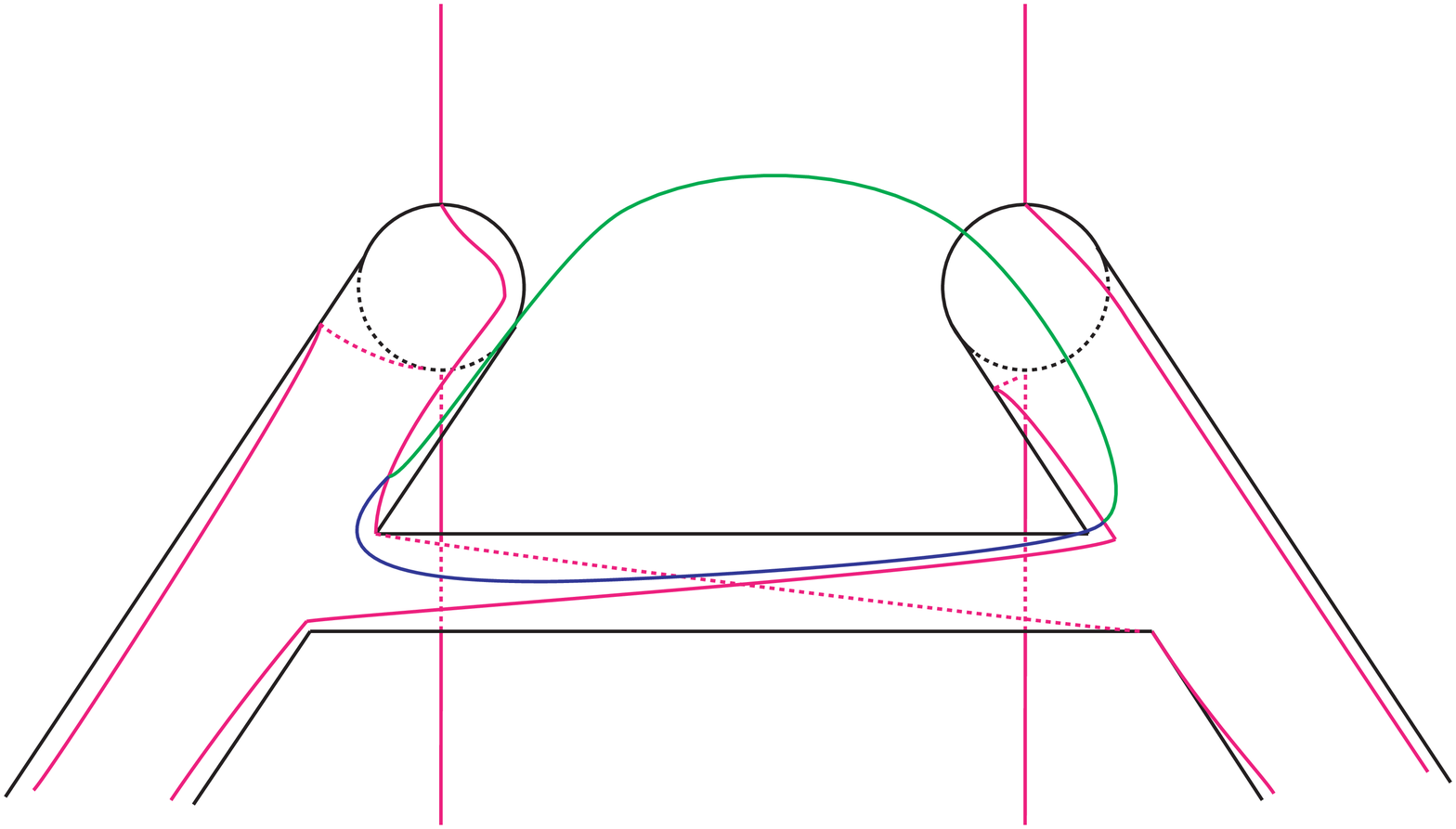}
\put(10,45){\tiny $F_1=\bdry M$} \put(71.5,52){\tiny $\Gamma_{\bdry
M}$} \put(50,46.5) {\tiny $d_i$} \put(20.7,19){\tiny $c_i$}
\end{overpic}
\caption{Step (1) of the subdivision process. Attaching an arc of
the first type. The face $F_1=\bdry M$ is in the back. The cylinders
on the left and right are thickenings of $\{p\}\times[0,1]$ and
$\{q\}\times[0,1]$, and the horizontal cylinder is a thickening of
$c_i$. The blue arc is $c_i$ and the green arc is $d_i$.}
\label{attacharc1}
\end{figure}
In this case there is an arc $d_i$ which passes through
$R_\pm(\Gamma)$ and is isotopic to $c_i$ rel endpoints inside
$M-N(K_i)$, so that $\gamma_i=c_i\cup d_i$ bounds a disk in
$M-N(K_i)$ and has $tb(\gamma_i)=-1$ with respect to this disk.
Another way of thinking about the attachment of $c_i$ is to slide
both endpoints of $c_i$ towards $F_1=\bdry M$, along
$\Gamma_{\bdry(M-N(K_i))}$, so that both endpoints of $c_i$ are then
placed on $\Gamma_{\bdry M}$. Figure~\ref{attacharc} depicts
\begin{figure}[ht]
\begin{overpic}[width=6cm]{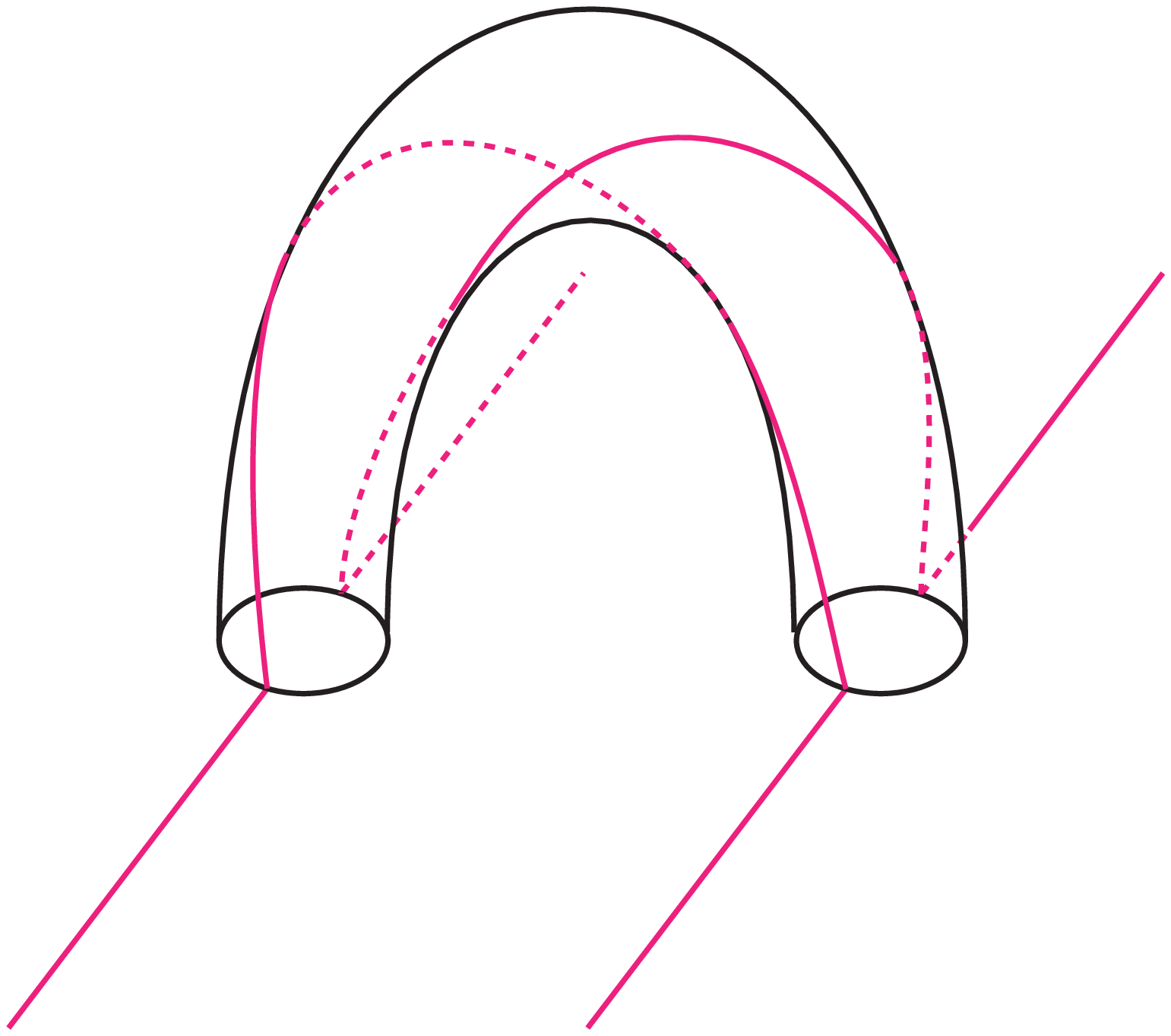}
\end{overpic}
\caption{The face $F_1=\bdry M$ is in the back.} \label{attacharc}
\end{figure}
the neighborhood of $c_i$ after the isotopy. The disk $D$ with
$tb(\bdry D)=-1$ is easy to see in this diagram.  Once the arcs of
the first type are attached, we need to attach arcs $c_i\subset
F_{1-\delta}$ of the second type in order to complete the Legendrian
$1$-skeleton of $F_{1-\delta}$. These $c_i$ connect interior points
of arcs of the first type and do not intersect
$\Gamma_F\times\{1-\delta\}$. This is given in
Figure~\ref{attacharc3}.  The top endpoint of $c_i$ in
Figure~\ref{attacharc3} can be moved to the left and the bottom
endpoint can be moved to the right, both along $\Gamma_{\bdry
(M-N(K_i))}$, so that both endpoints now lie on $\Gamma_{\bdry M}$.
The result is the same situation as given in Figure~\ref{attacharc}.
\begin{figure}[ht]
\begin{overpic}[width=6cm]{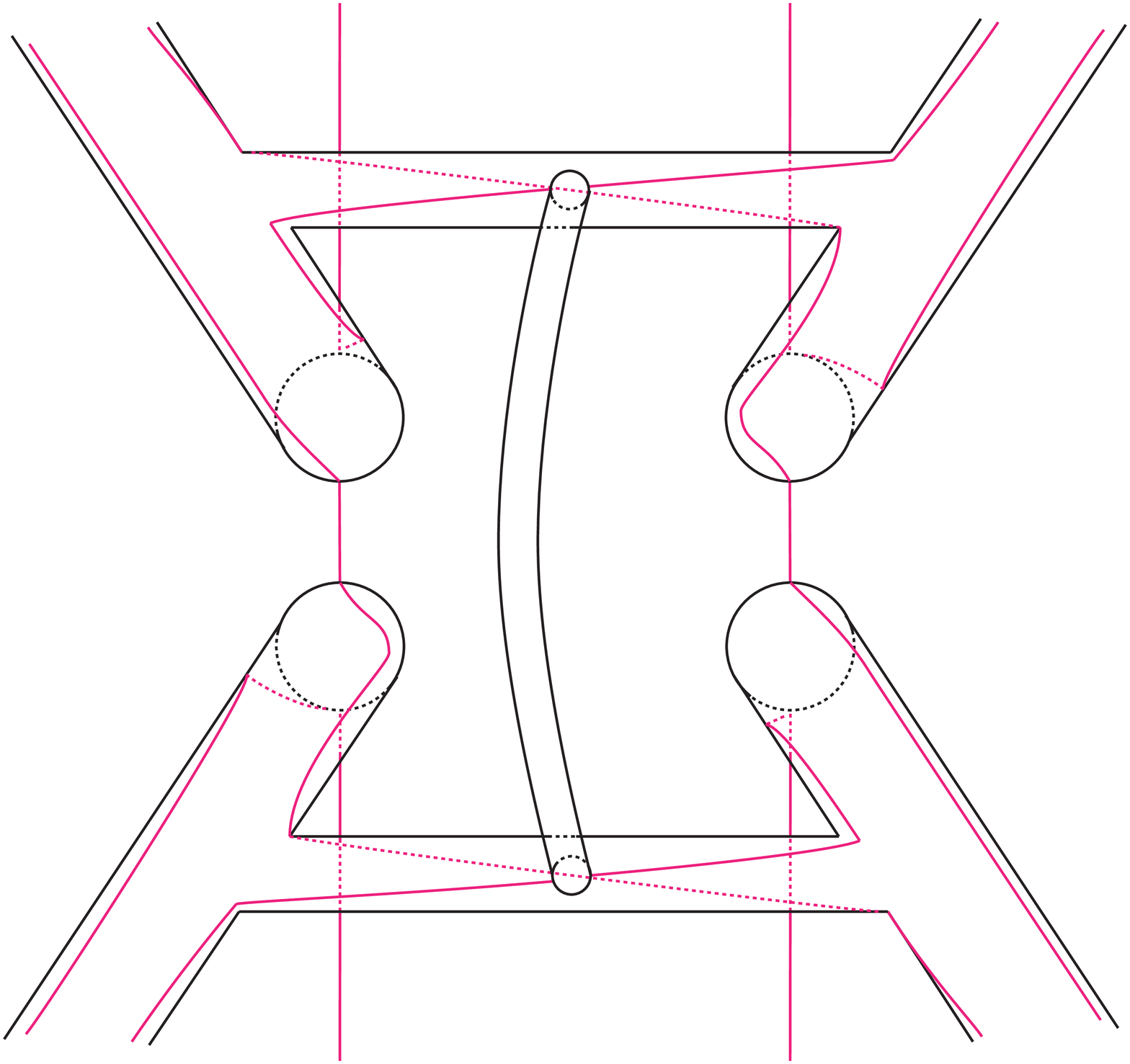}
\put(50,48){\tiny $N(c_i)$}
\end{overpic}
\caption{Step (1) of the subdivision process. Attaching an arc of
the second type. The face $F_1=\bdry M$ is in the back.}
\label{attacharc3}
\end{figure}

Next we consider Step (3), which is depicted in
Figure~\ref{attacharc2}. In the figure, the face $F_1=\bdry M$ is in
the back, and the thickening of the Legendrian $1$-skeleton of
$F_{1-\delta}$ is to the front. The thickening of the arc $c_i$ is
the cylinder to the left emanating from $\Gamma_{\bdry M}$.
\end{proof}

\section{Definition of the contact class} \label{section: defn}

We briefly recall Juh\'asz' sutured Floer homology theory
\cite{Ju1,Ju2}.

A {\em compatible} Heegaard splitting for a sutured manifold
$(M,\Gamma)$ consists of a Heegaard surface $\Sigma$ (not
necessarily connected) with nonempty boundary, together with two
sets of pairwise disjoint simple closed curves that do not intersect
$\bdry \Sigma$, the $\alpha$-curves $\alpha_1,\dots,\alpha_r$ and
the $\beta$-curves $\beta_1,\dots,\beta_r$.   Then $M$ is obtained
from $\Sigma\times[-1,1]$ by gluing compressing disks along
$\alpha_i\times \{-1\}$ and along $\beta_i\times\{1\}$, and
thickening.  We take the suture $\Gamma$ to be
$\bdry\Sigma\times\{0\}$.

Let $\mathbb{T}_\alpha= \alpha_1\times\dots\times \alpha_r$ and
$\mathbb{T}_{\beta}= \beta_1\times\dots\times \beta_r$, viewed in
$Sym^r(\Sigma)$. Then let $\CF(\Sigma,\alpha,\beta)$ be the free
$\Z$-module generated by the points $\mathbf{x}=(x_1,\dots,x_r)$ in
$\mathbb{T}_\alpha\cap \mathbb{T}_\beta$.  The suture $\Gamma$ plays
the role of the basepoint in sutured Floer homology.  Denote by
$\mathcal{M}_{\mathbf{x} ,\mathbf{y}}$ the 0-dimensional (after
quotienting by the natural $\R$-action) moduli space of holomorphic
maps $u$ from the unit disk $D^2\subset \C$ to $Sym^r(\Sigma)$ that
(i) send $1\mapsto \mathbf{x}$, $-1\mapsto \mathbf{y}$, $S^1 \cap
\{\mbox{Im } z \geq 0 \}$ to $\mathbb{T}_{\alpha}$ and $S^1 \cap
\{\mbox{Im } z \leq 0\}$ to $\mathbb{T}_{\beta}$, and (ii) avoid
$\bdry\Sigma \times Sym^{r-1}(\Sigma)\subset Sym^r(\Sigma)$. Then
define
$$\bdry \mathbf{x} =  \sum_{\mu(\mathbf{x},\mathbf{y})=1} ~~
\#(\mathcal{M}_{\mathbf{x}, \mathbf{y} })~~  \mathbf{y},$$ where
$\mu(\mathbf{x},\mathbf{y})$ is the relative Maslov index of the
pair and $\#(\mathcal{M}_{\mathbf{x}, \mathbf{y} })$ is a signed
count of points in $\mathcal{M}_{\mathbf{x}, \mathbf{y} }$. The
homology $\SFH(M,\Gamma)$ of this complex is shown to be independent
of the various choices made in the definition. In particular, it is
independent of the choice of a ``weakly admissible'' Heegaard
decomposition.

\s\n {\bf Interlude on orientation conventions.} The convention for
$(M,\Gamma)$ (consistent with that of Ozsv\'ath-Szab\'o and
Juh\'asz) is as follows:  The Heegaard surface $\Sigma$ is an
oriented surface whose oriented boundary is $\Gamma$. The
suture/dividing set $\Gamma$ is also defined as the boundary of
$R_+(\Gamma)$ ($=$ minus the boundary of $R_-(\Gamma)$). If $\Sigma$
splits $M$ into two compression bodies $H_1$ and $H_2$, and $\bdry
H_1=\Sigma$, $\bdry H_2=-\Sigma$, then the boundaries of the
compressing disks for $H_1$ are the $\alpha_i$ and the boundaries of
the compressing disks for $H_2$ are the $\beta_i$. Then
$\SFH(M,\Gamma)$ is the Floer homology $\SFH(\Sigma,\alpha,\beta)$,
in that order. We will now describe $\SFH(M,-\Gamma)$,
$\SFH(-M,\Gamma)$, and $\SFH(-M,-\Gamma)$, using the same data
$(\Sigma,\alpha,\beta)$. The reader is warned that the
$\Sigma,\alpha,\beta$ that appear in this interlude will be
different from the $\Sigma,\alpha,\beta$ that appear subsequently.

$\SFH(M,-\Gamma)$:  If we keep the same orientation for $M$ and
switch the orientation of $\Sigma$, then we must switch $\alpha$ and
$\beta$. Hence $\SFH(M,-\Gamma)= \SFH(-\Sigma,\beta,\alpha)$.  Also
$R_\pm(\Gamma)=R_\mp(-\Gamma)$.

$\SFH(-M,\Gamma)$:  If we switch the orientation of $M$, then $\bdry
H_1=-\Sigma$ and $\bdry H_2=\Sigma$.  Since the orientation of
$\Gamma$ is unchanged, the orientation of $\Sigma$ is unchanged.
Hence $\SFH(-M,\Gamma)=\SFH(\Sigma,\beta,\alpha)$.  Observe that
$R_\pm(M,\Gamma)= R_\mp(-M,\Gamma)$.

$\SFH(-M,-\Gamma)= \SFH(-\Sigma,\alpha,\beta)$, from the above
considerations.

\s Given $(M,\Gamma,\xi)$, the decomposition of $M$ into $M-N(K)$
and $N(K)$ from Theorem~\ref{thm: decomposition} gives us a {\em
partial open book decomposition} $(S,R_+(\Gamma),h)$, which we now
describe. (The terminology {\em partial} open book decomposition is
used because $M$ can be constructed from a page $S$ and a
partially-defined monodromy map
$h:P=S-\overline{R_+(\Gamma)}\rightarrow S$.)  The tubular portion
$T$ of $-\bdry N(K)$ is split by the dividing set into positive and
negative regions, with respect to the orientation of $-\bdry N(K)$
or $\bdry (M-N(K))$.  Let $P$ be the positive region. Next, if
$D_i\subset \bdry N(K)$ are the attaching disks of $N(K)$, then
consider $R_+(\Gamma)-\cup_i D_i$; from now on this subsurface of
$\bdry M$ will be called $R_+(\Gamma)$.  Then the page $S$ is
obtained from the closure $\overline{R_+(\Gamma)}$ of $R_+(\Gamma)$
by attaching the positive region $P$.    Let $\sim$ be an
equivalence relation on $S\times[-1,1]$ given by $(x,t) \sim
(x,t')$, where $x \in \bdry S$ and $t,t' \in [-1,1]$. Then
$S\times[-1,1]/_\sim$ can be identified with $M-N(K)$, where
$\Gamma_{\bdry (M-N(K))}=\bdry S\times\{0\}$. The manifold $M$ is
then obtained from $S\times[-1,1]/_\sim$ by attaching thickenings of
the compressing disks $D_i^\beta$ corresponding to the meridians of
$N(K)$. The suture $\Gamma$ on $\bdry M$ is $\bdry
(\overline{R_+(\Gamma)})\times\{1\}$ --- we emphasize that this
$\Gamma$ is no longer quite the same as the dividing set
$\Gamma_{\bdry M}$. Also observe that $\bdry D_i^\beta \cap (\bdry
(\overline{R_+(\Gamma)})\times\{1\})=\emptyset$. The Heegaard
surface is $\Sigma= \bdry (S\times[-1,1]/_{\sim}) -
(R_+(\Gamma)\times \{1\})$.  One set of compressing disks is
$\{D_i^\beta\}$ and the other set $\{D_i^\alpha\}$ gives a disk
decomposition of $M-N(K)$. Let $h$ be the monodromy map
--- it is obtained by first pushing $P$
across $N(K)$ to $T-P \subset \bdry (M-N(K))$, and then following it
with an identification of $M-N(K)$ with $S\times[-1,1]/_\sim$.

\begin{figure}[ht]
\begin{overpic}[width=12cm]{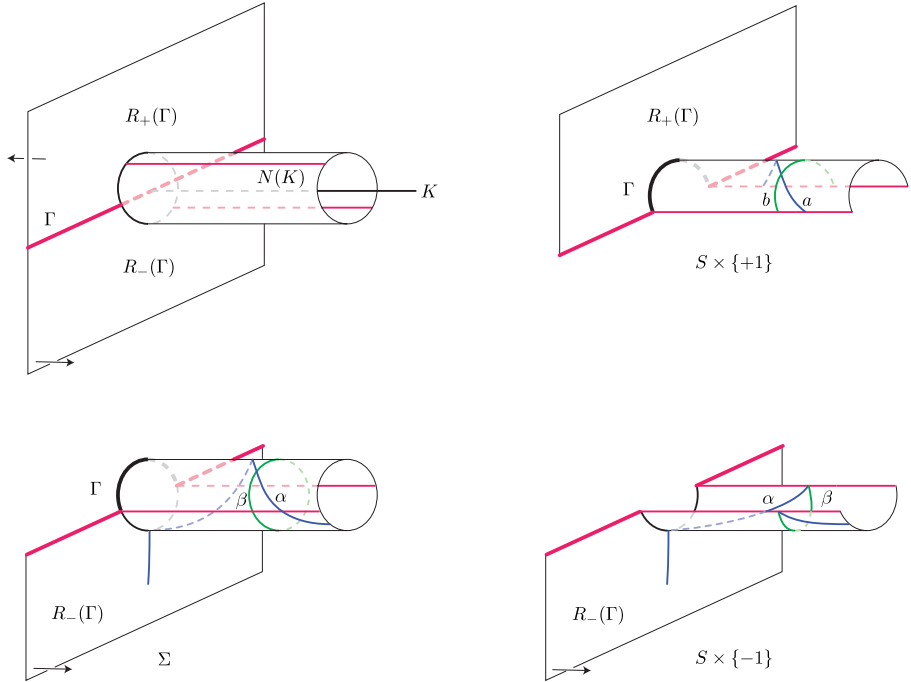}
\end{overpic}
\caption{The first figure shows a small portion of $M$ near an
intersection point of $\Gamma$ and $K$ and is drawn from the point
of view of the interior of $M-N(K)$. The second figure shows the
page $S\times\{1\}$ and a basis arc $a$. Notice that in this figure
$R_+(\Gamma)$ denotes a proper subsurface of the previous
$R_+(\Gamma)$.  The boundary of the new $R_+(\Gamma)$ is the new
$\Gamma$.  We have also applied edge-rounding to connect the
dividing curves on $\bdry M-\cup_i D_i$ to those on the tubular
portion $T$. The third figure shows the Heegaard surface $\Sigma$.
Notice that the arc $\beta - b$ is obtained by pushing $b$ through
$N(K)$, while $\alpha-a$ is obtained by isotoping $a$ through
$M-N(K)$. The subsurface of $\Sigma$ indicated in the fourth figure
is $S\times\{-1\}$. Notice that the single point of intersection
$a\cap b$ is replaced by two copies of itself in this figure. In
doing sutured Floer homology computations, we only consider
holomorphic disks which miss complementary regions of $\alpha \cup
\beta$ containing $\Gamma$.  This ensures that the computations may
be done in the subsurface $S\times\{-1\}$ of $\Sigma$.}
\label{manifold}
\end{figure}

\begin{figure}[ht]
\begin{overpic}[width=12cm]{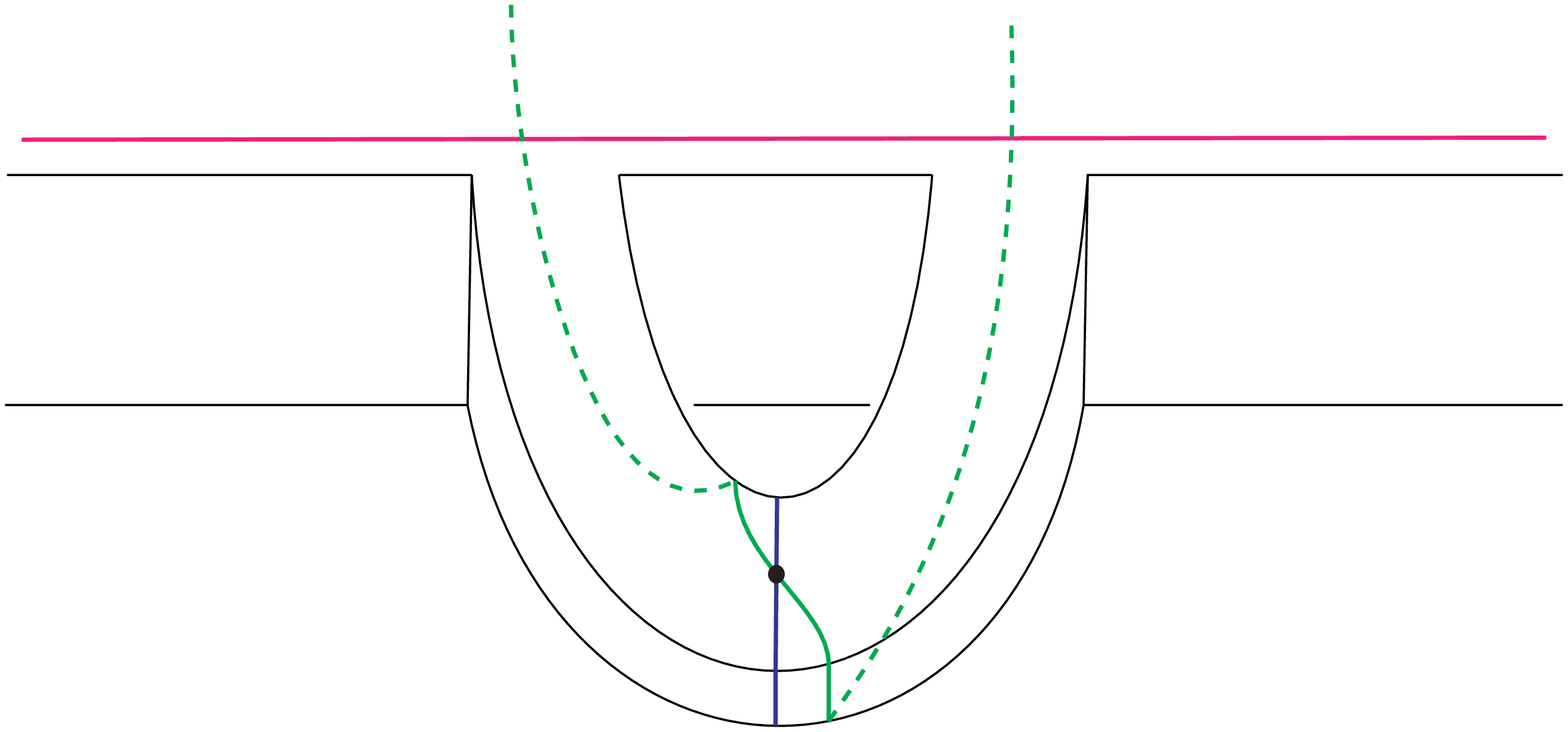}
\put(10,38){\tiny $\Gamma$} \put(51.2,9.3){\tiny $x$}
\put(47.5,6){\tiny $\alpha$} \put(59.8,18.2){\tiny $\beta$}
\end{overpic}
\caption{The partial open book decomposition. The curve $\Gamma$ is
the suture on $\bdry M$. (It has been pushed into $R_+(\Gamma)$ for
better viewing.) The top face is $S\times\{1\}$ and contains
$R_+(\Gamma)$. In this figure, a single 1-handle is attached to
$\overline{R_+(\Gamma)}$. The point $x$ is then the contact class.}
\label{page}
\end{figure}

Suppose that $S$ is obtained by successively attaching $r$ 1-handles
to the union of $\overline{R_+(\Gamma)}$ and the previously attached
1-handles. Let $a_i$, $i=1,\dots,r$, be properly embedded arcs in
$P=S-\overline{R_+(\Gamma)}$ with endpoints on $A=\bdry P -\Gamma$,
so that $S-\cup_i a_i$ deformation retracts onto
$\overline{R_+(\Gamma)}$. A collection $\{a_1,\dots,a_r\}$ of such
arcs is called a {\em basis} for $(S,R_+(\Gamma))$. In fact,
$\{a_1,\dots,a_r\}$ is a basis for $H_1(P,A)$. Next let $b_i$ be an
arc which is isotopic to $a_i$ by a small isotopy so that the
following hold: \be
\item The endpoints of $a_i$ are isotoped along $\bdry S$, in the
direction given by the boundary orientation of $S$.
\item The arcs $a_i$ and $b_i$ intersect transversely in one point in the
interior of $S$.
\item If we orient $a_i$, and $b_i$ is given the induced orientation
from the isotopy, then the sign of the intersection $a_i\cap b_i$ is
$+1$. \ee Then the $\alpha$-curves are $\bdry (a_i\times[-1,1])$ and
the $\beta$-curves are $(b_i\times\{1\})\cup (h(b_i)\times\{-1\}),$
viewed on $S\times[-1,1]/_\sim$. The $\alpha$-curves and
$\beta$-curves avoid the suture $\Gamma$ and hence determine a
Heegaard splitting $(\Sigma,\beta,\alpha)$, which is easily seen to
be weakly admissible. Now the union of $\Sigma$ and
$R_+(\Gamma)=R_+(\Gamma)\times\{1\}$ bound $S\times[0,1]/\sim$.
Define the orientation on $\Sigma$ to be the outward orientation
inherited from $S\times[0,1]/\sim$.  Then $\bdry \Sigma$ is oriented
oppositely from $\bdry(R_+(\Gamma))$, and we have
$\bdry\Sigma=-\Gamma$. Hence $SFH(\Sigma,\beta,\alpha)=
SFH(-M,-\Gamma)$.

The contact class is basically the $\EH$ class which was defined in
\cite{HKM2}. The only difference is that we are not using a full
basis for $S$ and that the contact class sits in
$\SFH(-M,-\Gamma)/\{\pm 1\}$. Let $\mathbb{T}_\alpha=
\alpha_1\times\dots\times \alpha_r$ and $\mathbb{T}_{\beta}=
\beta_1\times\dots\times \beta_r$, viewed in $Sym^r(\Sigma)$. Let
$\CF(\Sigma,\beta,\alpha)$ be the chain group generated by the
points in $\mathbb{T}_\beta\cap \mathbb{T}_\alpha$. Let $x_i$ be the
intersection point $(a_i\cap b_i)\times\{1\}$ lying in
$S\times\{1\}$.  Then $\mathbf{x}=(x_1,\dots,x_r)$ is a cycle in
$\CF(\Sigma,\beta,\alpha)$ due to the placement of $\Gamma$, and its
class in $\SFH(-M,-\Gamma)$ will be written as
$\EH(S,h,\{a_1,\dots,a_r\})$. Figures~\ref{manifold} and \ref{page}
depict the situation described above.

\section{Well-definition of the contact class} \label{section:
well-defined}

\begin{thm}
$\EH(S,h,\{a_1,\dots,a_r\})\in \SFH(-M,-\Gamma)/\{\pm 1\}$ is an
invariant of the contact structure.
\end{thm}

Once the invariance is established, the contact class will be
written as $\EH(M,\Gamma,\xi)$.   As usual, there is a $\pm 1$
indeterminacy, which we sometimes suppress.

\begin{proof}[Outline of Proof.]
We need to prove that $\EH(S,h,\{a_1,\dots,a_r\})$ is (i) independent
of the choice of basis $\{a_1,\dots,a_r\}$, (ii) invariant under
stabilization, and (iii) only depends on the isotopy class of $h$.
The dependence only on the isotopy class of $h$ is identical to the
proof given in \cite{HKM2}, and will be omitted.  The proof of the
independence of choice of basis is similar to that of \cite{HKM2},
and we highlight only the differences in Section~\ref{subsection:
basis}.

We define a {\em positive stabilization} in the relative case as
follows: Let $c$ be a properly embedded arc in $S$; in particular,
it can pass through $R_+(\Gamma)$. Attach a $1$-handle to $S$ at the
endpoints of $c$ to obtain $S'$, and let $\gamma$ be a closed curve
obtained by gluing $c$ and a core of the $1$-handle. Then a {\em
positive stabilization} of $(S,R_+(\Gamma),h)$ is $(S',R_+(\Gamma),
R_\gamma\circ h)$, where $R_\gamma$ is a positive Dehn twist about
$\gamma$.  We emphasize that the order of composing $R_\gamma$ and
$h$ is important, since we are not allowed a global conjugation of
$S$ when we only have a partial open book.

Every two partial open book decompositions representing the same
contact structure $\xi$ become isotopic after performing a sequence
of positive stabilizations to each. The proof follows from
Theorem~\ref{thm: subdivision}: Let us consider the case when $c_0$
is attached to $K_0$.  Since $M-N(K)\simeq S\times[-1,1]/_\sim$,
Condition~\ref{unknot} of Theorem~\ref{thm: subdivision} implies
that $c_0$ can be viewed as a Legendrian arc on $S\times\{0\}\subset
S\times[-1,1]/_\sim$. Attaching a neighborhood of $c_0$ to $N(K_0)$
is equivalent to drilling out a standard neighborhood of the
Legendrian arc $c_0\subset S\times\{0\}$. Recall that the monodromy
map $h$ sends $P\times\{1\}$ to $S\times\{-1\}$. Since the drilling
takes place in the region $S\times[-1,1]$, {\em after} $h$ is
applied, the new monodromy map is $R_\gamma\circ h$.

The invariance under stabilization in the closed case can be argued
as follows, once we establish the independence of choice of basis:
Suppose the stabilization occurs along the properly embedded arc $c$
in $S$.  Take a basis $\{a_1,\dots,a_r\}$ for $S$ so that all the
$a_i$ are disjoint from $c$.  The arc $c$ may be separating or
nonseparating (even boundary-parallel), but there is a suitable
basis in either case. Then take the basis $\{a_0,a_1,\dots,a_r\}$
for $S'$, where $a_0$ is the cocore of the 1-handle attached onto
$S$. Then $\beta_0$ and $\alpha_0$ intersect exactly once, and
$\beta_0$ is disjoint from the other $\alpha_i$ by construction.
Although there is a natural chain isomorphism between
$\CF(\beta,\alpha)$ and
$\CF(\beta\cup\{\beta_0\},\alpha\cup\{\alpha_0\})$, where
$\alpha=\{\alpha_1,\dots,\alpha_r\}$ and
$\beta=\{\beta_1,\dots,\beta_r\}$, it is not clear that the
identification is consistent with the stabilization and handleslide
maps in \cite{OS1}. In Lemma~\ref{lemma: complexity zero} we prove
that $\EH(S,h,\{a_1,\dots,a_r\})$ is indeed mapped to
$\EH(S',h',\{a_0,\dots,a_r\})$ under the stabilization and
handleslide maps in \cite{OS1}. (Notice that this proof of
invariance could have been used in \cite{HKM2} without appealing to
the equivalence with the Ozsv\'ath-Szab\'o contact invariant. This
was pointed out to the authors by Andr\'as Stipsicz.)

Now, if $R_+(\Gamma)$ is not a homotopically trivial disk, then it
is not always possible to find a basis $\{a_1,\dots,a_r\}$ which is
disjoint from the arc of stabilization $c$.  This difficulty is
dealt with in Section~\ref{subsection: stab}.
\end{proof}

\subsection{Change of basis}\label{subsection: basis}

Let $\{a_1,a_2,\dots,a_r\}$ be a basis for $(S,R_+(\Gamma))$. After
possibly reordering the $a_i$'s, suppose $a_1$ and $a_2$ are
adjacent arcs on $A=\bdry P-\Gamma$, i.e., there is an arc
$\tau\subset A$ with endpoints on $a_1$ and $a_2$ such that $\tau$
does not intersect any $a_i$ in $\mbox{int}(\tau)$. Define $a_1+a_2$
as the isotopy class of $a_1\cup\tau\cup a_2$, relative to the
endpoints. Then the modification $\{a_1,a_2,\dots,a_r\}\mapsto
\{a_1+a_2,a_2,\dots,a_r\}$ is called an {\em arc slide}.

\begin{lemma}  \label{independence}
$\EH(S,h)$ is invariant under an arc slide
$\{a_1,a_2,\dots,a_r\}\mapsto \{a_1+a_2,a_2,\dots,a_r\}$.
\end{lemma}

\begin{proof}
Same as that of Lemma~3.4 of \cite{HKM2}.
\end{proof}

Let $\{a_1,\dots,a_r\}$ and $\{b_1,\dots,b_r\}$ be two bases for
$(S,R_+(\Gamma))$. Assume that the two bases intersect transversely
and {\em efficiently}, i.e., each pair of arcs $a_i, b_j$ realizes
the minimum number of intersections in its isotopy class, where the
endpoints of the arcs are allowed to move in $A$. In particular,
there are no bigons consisting of a subarc of $a_i$ and a subarc of
$b_j$, and no triangles consisting of a subarc of $a_i$, a subarc of
$b_j$, and a subarc of $A$.

\begin{lemma}
Suppose that each component of $\overline{P}$ intersects $\Gamma$
along at least two arcs. Then there is a sequence of arc slides
which takes $\{a_1,\dots,a_r\}$ to $\{b_1,\dots,b_r\}$.
\end{lemma}

Here $\overline{P}$ is the closure of $P$. The condition of the
lemma is easily satisfied by performing a trivial stabilization
along a boundary-parallel arc. Moreover, a trivial stabilization is
easily seen to preserve the $\EH$ class.

\begin{proof}
Consider a connected component $Q$ of $S- \cup_{i=1}^r a_i-
\overline{R_+(\Gamma)}$. Then $Q$ is a (partially open) polygon
whose boundary $\bdry Q$ consists of $2k$ arcs, $k-1$ of which are
$a_i$ or $a_i^{-1}$, $k$ of which are subarcs
$\tau_1,\dots,\tau_{k}$ of $\bdry P-\Gamma$, and one which is a
subarc $\gamma$ of $\Gamma$. (For the moment we have oriented the
$a_i$, and the notation $a_i^{-1}$ means that the orientation of
$a_i$ and the orientation of $\bdry Q$ are opposite.) Since each
component of $\overline{P}$ intersects $\Gamma$ along at least two
arcs, there is at least one arc $c\subset \bdry Q$ of type $a_i$ or
$a_i^{-1}$ so that $c^{-1}$ does not appear on $\bdry Q$. Otherwise,
$Q$ glues up to a component of $P$ whose closure intersects $\Gamma$
along one arc.

Suppose $(\bigcup_{i=1}^r a_i)\cap (\bigcup_{i=1}^r b_i)\not=0$.  We
will apply a sequence of arc slides to $\{a_i\}_{i=1}^r$ to obtain
$\{a_i'\}_{i=1}^r$ which is disjoint from $\{b_i\}_{i=1}^r$. After
possibly reordering the arcs, there is a subarc $b_1^0\subset b_1$
in $Q$ with endpoints on $a_1$ and $\tau_1$ on $\bdry Q$.  The
subarc $b_1^0$ separates $Q$ into two regions $Q_1$ and $Q_2$, only
one of which (say $Q_1$) has $\gamma$ as a boundary arc.  If
$c=a_1$, then we can slide $c$ around the boundary of $Q_2$ to
obtain $a_1'$ which has fewer intersections with $b_1$. If $a_1$ and
$a_1^{-1}$ both occur on $\bdry Q$, then we have a problem if $Q_1$
has $\gamma$ as a boundary arc and $Q_2$ has $a_1^{-1}$ as a
boundary arc.  To maneuver around this problem, move $c$ via a
sequence of arc slides so that the resulting $c'$ is parallel to
$\gamma$ (and protects it). Then we can slide $a_1$ around the
boundary of $\bdry Q_1$ as in \cite{HKM2}.

Finally, suppose that $\{a_i\}_{i=1}^r$ and $\{b_i\}_{i=1}^r$ are
disjoint. Let $Q$ be a component of $S-\cup_i
a_i-\overline{R_+(\Gamma)}$, as before.  Let $b_1$ be an arc in $Q$
which is not parallel to any $a_i$ or $a_i^{-1}$.  The arc $b_1$
cuts $Q$ into two components $Q_1$ and $Q_2$, where $Q_1$ has
$\gamma$ as a boundary arc.  We claim that there is some $a_i$ on
$\bdry Q_2$ so that $a_i^{-1}$ is not on $\bdry Q_2$.  Otherwise,
the arcs on $\bdry Q_2$ will be paired up and $b_1$ will cut off a
subsurface of $S$ whose closure does not intersect $\Gamma$.
Moreover, we may also assume that $a_i$ is not parallel to any
$b_j$.  Now, apply a sequence of arc slides to $a_i$ so that it
becomes parallel to $b_1$.
\end{proof}

\subsection{Stabilization} \label{subsection: stab}

A properly embedded arc $c$ in $S$ that intersects $\Gamma$
efficiently is said to have {\em complexity $n$} if there are $n$
subarcs of $c$ in $\overline{P}$, both of whose endpoints are on
$\Gamma$.  Let $(S',h')$ be a positive stabilization of a partial
open book decomposition $(S,h)$ along a properly embedded arc
$c\subset S$ that intersects $\Gamma$ efficiently. Then $(S',h')$ is
a stabilization of $(S,h)$ of {\em complexity $n$} if $c$ has
complexity $n$.

The following lemma gives the fundamental property of arcs of
complexity zero:

\begin{lemma} \label{lemma: 0}
An arc $c$ has complexity zero if and only if there is a basis for
$(S,R_+(\Gamma))$ which is disjoint from $c$.
\end{lemma}

\begin{proof}
If an arc $c$ has complexity $>0$, then there is a subarc
$c_0\subset c$ in $\overline{P}$, both of whose endpoints are on
$\Gamma$. It is impossible to find a basis that does not intersect
$c_0$, since the pair $(c_0,\bdry c_0)$ is homotopically nontrivial
in $(\overline{P},\Gamma)$.

On the other hand, suppose the complexity is zero. Then there are at
most two subarcs of $c\cap \overline{P}$. If $c$ lies in
$\overline{R_+(\Gamma)}$, then it does not intersect any basis of
$(S,R_+(\Gamma))$.  If $c$ is entirely contained in $P$, then it is
easy to see that there is a basis which avoids $c$. (There are two
cases, depending on whether $c$ cuts off a subsurface of $P$ whose
boundary does not intersect $\Gamma$.) Let $c_i$ be a subarc of
$c\cap \overline{P}$. Depending on $c$, there may be two such, i.e.,
$i=1,2$, or only one such, i.e., $i=1$. In either case, for each
$c_i$, one of its endpoints is on $\Gamma$ and the other on $\bdry
P-\Gamma$. If $c_i$ is a trivial subarc, i.e., forms a triangle in
$P$, together with an arc of $\Gamma$ and an arc of $\bdry
P-\Gamma$, then $c_i$ does not obstruct the formation of a basis and
can in fact be isotoped out of $P$. Therefore assume that $c_i$ is
nontrivial. Let $\tau_1, \gamma, \tau_2$ be subarcs of $\bdry P$ in
counterclockwise order, where $\gamma\subset \Gamma$ contains an
endpoint of $c_1$, and $\tau_i$ are subarcs of $\bdry P-\Gamma$.
Also let $\tau$ be the component of $\bdry P-\Gamma$ with the other
endpoint of $c_1$. Then take the first arc $a$ of the basis to be
parallel to $c_1$, with endpoints on either $\tau$ and $\tau_1$ or
$\tau$ and $\tau_2$. We have a choice of either of these unless
$c_2$ also has an endpoint on $\gamma$, in which case only one of
these will work (and not intersect $c_2$). The nontriviality of
$c_1$ implies that $a$ is neither boundary-parallel, i.e., cuts off
a disk in $P$, nor parallel to $\Gamma$. If the arc $a$ cuts off a
subsurface of $P$ whose boundary does not intersect $\Gamma$, then
we must replace it with disjoint, nonseparating arcs
$a_1,\dots,a_{2g}$, where $g$ is the genus of the cut-off surface.
Once we cut $P$ along $a$ or $a_1,\dots,a_{2g}$, then $c_1$ becomes
trivial in the cut-open surface, and can be ignored. We apply the
same technique to $c_2$, if necessary, and the rest of the basis can
be found without difficulty.
\end{proof}

\begin{lemma}\label{lemma: complexity zero}
Let $(S',h')$ be a stabilization of $(S,h)$ along an arc $c$ of
complexity zero. Suppose $\{a_1,\dots,a_r\}$ is a basis for $S$
which is disjoint from $c$ and $\{a_0,a_1,\dots,a_r\}$ is a basis
for $S'$, where $a_0$ is the cocore of the 1-handle attached to
obtain $S'$.  Then $\EH(S,h,\{a_1,\dots,a_r\})$ is mapped to
$\EH(S',h',\{a_0,\dots,a_r\})$ under the stabilization and
handleslide maps defined in \cite{OS1}.
\end{lemma}

\begin{figure}[ht]
\begin{overpic}[width=11cm]{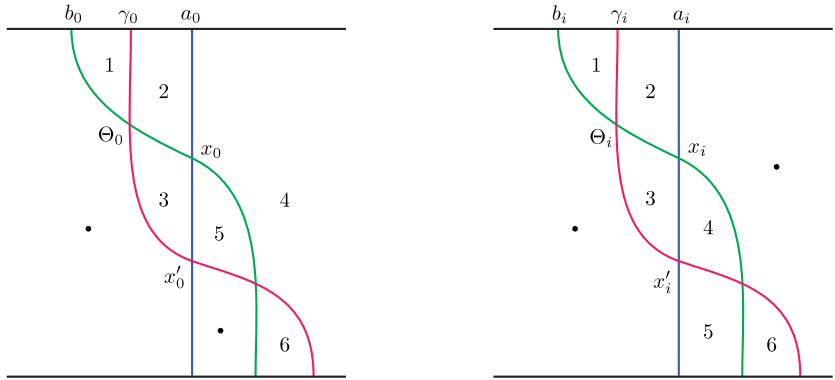}
\end{overpic}
\caption{The blue arc is $a_i$, the green arc is $b_i$, and the red
arc is a pushoff of $b_i$.} \label{tensoring}
\end{figure}

\begin{proof}
If no $\beta_i$, $i>0$, nontrivially intersects $\alpha_0$, we have
a standard stabilization. Suppose there is some $\beta_i$, $i>0$,
which nontrivially intersects $\alpha_0$.  This means $h(b_i)$
nontrivially intersected $c$ before stabilization. Let $\eta$ be a
closed curve on $S'$, obtained from $c$ by attaching the core of the
1-handle. In the left-hand diagram of Figure~\ref{handleslide}, a
portion of $h(b_i)$ intersecting $c$ is drawn; $h'(b_i)$ is obtained
from $h(b_i)$ by applying $R_\eta$. In order to remove its
intersection with $\alpha_0$, we handleslide $\beta_i$ across
$\beta_0$ to obtain $\gamma_i$, as indicated on the right-hand side
of Figure~\ref{handleslide}. In Figures~\ref{handleslide}
and~\ref{tensoring}, we place a dot in a region of $\bdry
(S\times[-1,1]/_\sim) -\cup_{i=0}^r \alpha_i-\cup_{i=0}^r \beta_i$
to indicate that any point in the region has a path in the region
that connects to $\Gamma=\bdry \Sigma$. In the right-hand diagram of
Figure~\ref{handleslide}, the dot to the left of $\alpha_0$ is
clear; by following $\beta_0$, the region with the dot to the right
of $\alpha_0$ is the same as the region with the dot to the left of
$\alpha_0$.

\begin{figure}[ht]
\begin{overpic}[width=16.5cm]{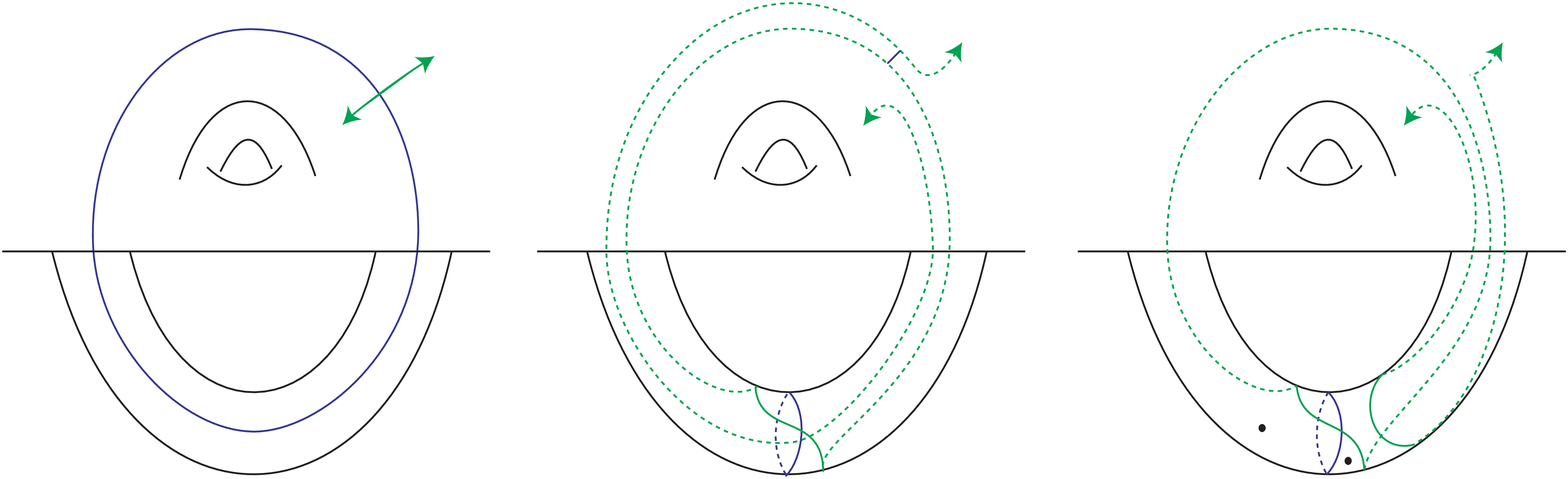}
\put(27,16){\tiny $\eta$} \put(27.5,25){\tiny $h(b_i)$}
\put(51.03,4.4){\tiny $\alpha_0$} \put(42.8,8.1){\tiny $\beta_0$}
\put(83.4,2){\tiny $x_0$}
\end{overpic}
\caption{The left-hand diagram shows the closed curve $\eta$ and an
intersection point of $h(b_i)$ and $\eta$, on a page $S$.  The
middle diagram is the Heegaard diagram before the handleslide, and
right-hand diagram gives the Heegaard diagram after the handleslide.
The dotted curves are on the page $S\times\{-1\}$ and the solid
curves are on $S\times\{1\}$. The short blue arc is on
$S\times\{-1\}$ and is the handleslide arc. } \label{handleslide}
\end{figure}
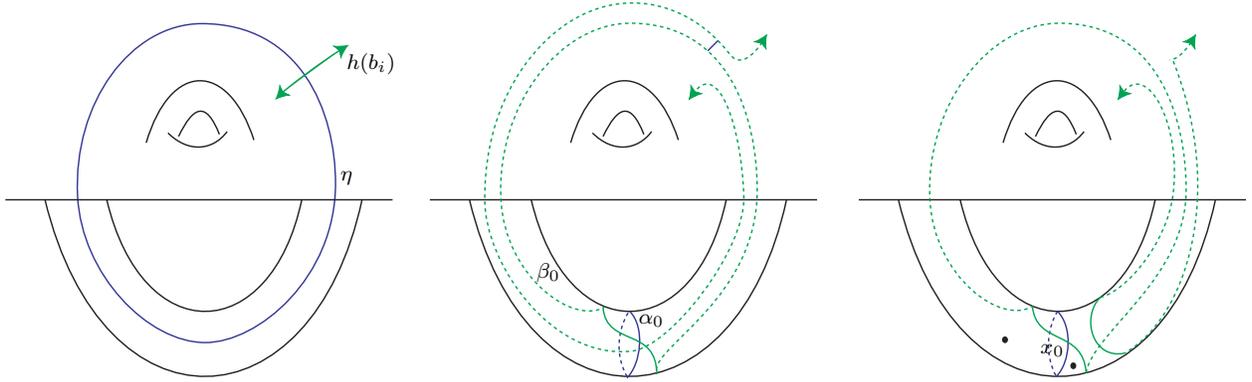

Let $\gamma_j$, $j\not=i$, be small pushoffs of $\beta_j$.
Figure~\ref{tensoring} shows the corresponding arcs in $P$. The
left-hand side of Figure~\ref{tensoring} corresponds to $i=0$ (and
Figure~\ref{handleslide}), and the right-hand side corresponds to
$i\not=0$.  We claim that the triple-diagram
$(\Sigma,\gamma,\beta,\alpha,\Gamma)$ is weakly admissible. Recall
that a triple-diagram is {\em weakly admissible} if each nontrivial
triply-periodic domain which can be written as a sum of
doubly-periodic domains has both positive and negative coefficients.
Consider the situation on the right-hand side of
Figure~\ref{tensoring}. The connected components of
$\Sigma-\cup_i\alpha_i-\cup_i\beta_i-\cup_i\gamma_i$ are numbered as
in the diagram. Due to the placement of $\Gamma$, the only potential
doubly-periodic domain involving $\beta,\alpha$ in the vicinity of
the right-hand diagram is $D_1+D_2-D_4-D_5$. Similarly, the only
domain involving $\gamma,\beta$ is $D_1+D_6-D_3-D_4$, and the only
one involving $\gamma,\alpha$ is $D_2+D_3-D_5-D_6$. Taking linear
combinations, we have
\begin{eqnarray*}
a(D_1+D_2-D_4-D_5)+b(D_1+D_6-D_3-D_4)+c(D_2+D_3-D_5-D_6) \\
= (a+b)D_1+ (a+c)D_2 +(-b+c) D_3 -(a+b)D_4-(a+c) D_5 +(b-c)D_6.
\end{eqnarray*}
Since the coefficients come in pairs, e.g., $a+b$ and $-(a+b)$, if
any of $a+b$, $b+c$, $a-c$ does not vanish, then the triply-periodic
domain has both positive and negative coefficients. Hence, if any of
$\alpha_i$, $\beta_i$ and $\gamma_i$ is used for some $i>0$, then we
are done. Otherwise, we may assume that none of $\alpha_i$,
$\beta_i$ and $\gamma_i$ is used in the periodic domain, for any
$i>0$. This allows us to assume that only $\alpha_0$, $\beta_0$ and
$\gamma_0$ are used for the boundary of the periodic domain. Once
all the $\alpha_i,\beta_i,\gamma_i$, $i>0$, are ignored, the $D_4$
region of the left-hand diagram becomes connected to $\Gamma$, and
we can apply the same procedure.

We now tensor $\mathbf{x}=(x_0,\dots,x_r)\in \CF(\beta,\alpha)$ with
the top generator $\Theta\in \CF(\gamma,\beta)$.  We claim that we
obtain $\mathbf{x'}=(x_0',\dots,x_r')$, as depicted in
Figure~\ref{tensoring}. As in the proof of the weak admissibility,
we start with the vicinity of right-hand side, i.e., $i>0$. There
are no other holomorphic triangles besides the obvious small
triangle involving $\Theta_i$, $x_i$, and $x_i'$, due to the
placement of the dots.  Once the $\alpha_i,\beta_i,\gamma_i$, $i>0$,
are exhausted, $D_4$ becomes connected to $\Gamma$, and we only have
the small triangle with vertices $\Theta_0$, $x_0$, $x_0'$.

Let the new $\beta$ be the old $\gamma$ and let the new $\mathbf{x}$
be the old $\mathbf{x'}$.  After a sequence of such handleslides
(and renamings), we obtain $\alpha$ and $\beta$ so that $\alpha_0$
and $\beta_0$ intersect once and do not intersect any other
$\alpha_i$ and $\beta_i$.  The standard stabilization map simply
forgets $x_0$.
\end{proof}

\begin{prop} \label{prop: stab}
Let $(S',h')$ be a stabilization of $(S,h)$ along an arc $c$ of
complexity $n$. Then there exists a stabilization of complexity at
most $n-1$ taking $(S',h')$ to $(S'',h'')$, and two stabilizations
of complexity at most $n-1$ which take $(S,h)$ to $(S'',h'')$.
\end{prop}

\begin{proof}
A stabilization corresponds to drilling out a neighborhood of a
properly embedded Legendrian arc $c=c\times\{0\}\subset
S\times\{0\}$ inside $S\times [-1,1]/_\sim$. We may assume that
$\Gamma_{\bdry (S\times[-1,1]/_\sim)}=\bdry S\times\{0\}$ and that
$S\times\{0\}$ is a convex surface with dividing set
$\Gamma_{S\times\{0\}}=\bdry S\times\{0\}$.  Then $c$ can be
realized as a Legendrian arc on $S\times\{0\}$ after applying the
{\em Legendrian realization principle} (see e.g. ~\cite{H1}). The
nonisolating condition, required in the Legendrian realization
principle, is easily met. The arc $c$ may be thought of having
twisting number $-{1\over 2}$ relative to $S\times\{0\}$, since it
begins and ends on the dividing set and intersects no other dividing
curves. Our strategy is to attach a Legendrian arc $d\subset
P\times\{0\}$ with one endpoint on $(\bdry P-\Gamma)\times \{0\}$
and the other endpoint on $c$ to obtain a Legendrian graph $L=c\cup
d$ on $S\times\{0\}$. (Again, the nonisolating condition is
satisfied.) The endpoint of $d$ on $c$ must lie on a component of
$c\cap \overline{P}$ with both endpoints on $\Gamma$. Then
$S''\times[-1,1]/_\sim$ will be $(S\times[-1,1]/_\sim) - N(L)$. We
will exploit two ways of obtaining the trivalent graph $L$.  One way
is to stabilize along $c$ and then along $d$.  For the other way,
subdivide $c=c_1\cup c_2$ at the trivalent vertex of $L$. We label
the arcs $c_1$, $c_2$ so that $d,c_1,c_2$ are in counterclockwise
order about the trivalent vertex. Then first stabilize along
$c_1\cup d$ and then along $c_2$. See Figure~\ref{graph}.

\begin{figure}[ht]
\begin{overpic}[width=8cm]{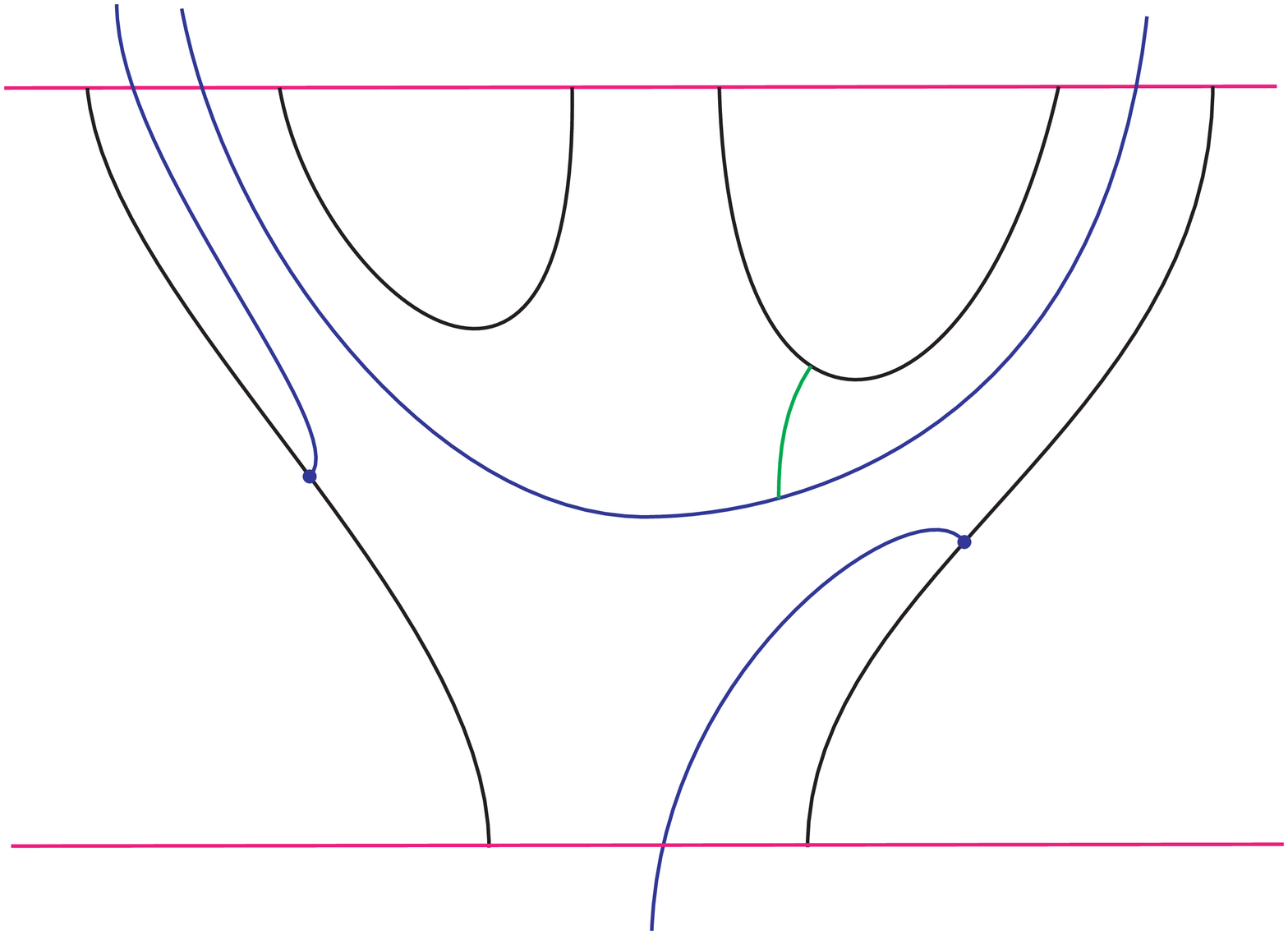}
\put(58,39){\tiny $d$} \put(45,30){\tiny $c_1$} \put(76,40){\tiny
$c_2$} \put(41,19){\tiny $P$} \put(-5,65.1){\tiny $\Gamma$}
\put(-5,6.3){\tiny $\Gamma$}
\end{overpic}
\caption{The red arcs are $\Gamma$ and the region drawn is $P$. }
\label{graph}
\end{figure}
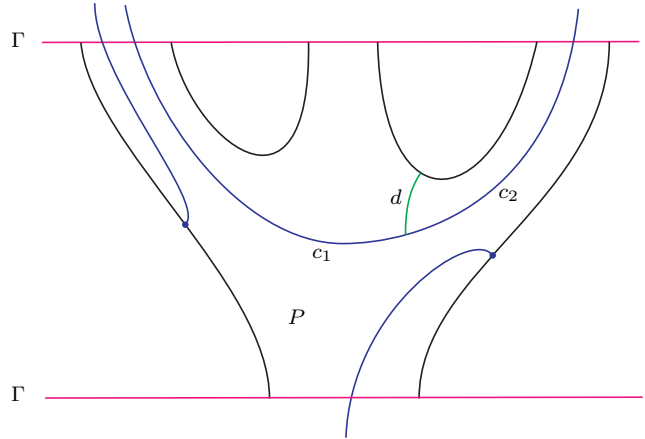

Next refer to Figure~\ref{stab}. We claim that if we stabilize along
the arc $c_1\cup d$ (its thickening is drawn in the diagram to the
left), and then stabilize along $c_2$, given as the green arc, then
the two stabilizations are equivalent to the stabilizations on the
right-hand side.  To see this, we isotop $c_2$ inside
$B=(S\times[-1,1]/_\sim) - N(c_1\cup d)$, while constraining its
endpoints to lie on $\Gamma_{\bdry B}$.  The key point to observe is
that $c_2$ can be slid ``to the right'' so that it lies above the
neighborhood of $c_1\cup d$.

\begin{figure}[ht]
\begin{overpic}[width=12cm]{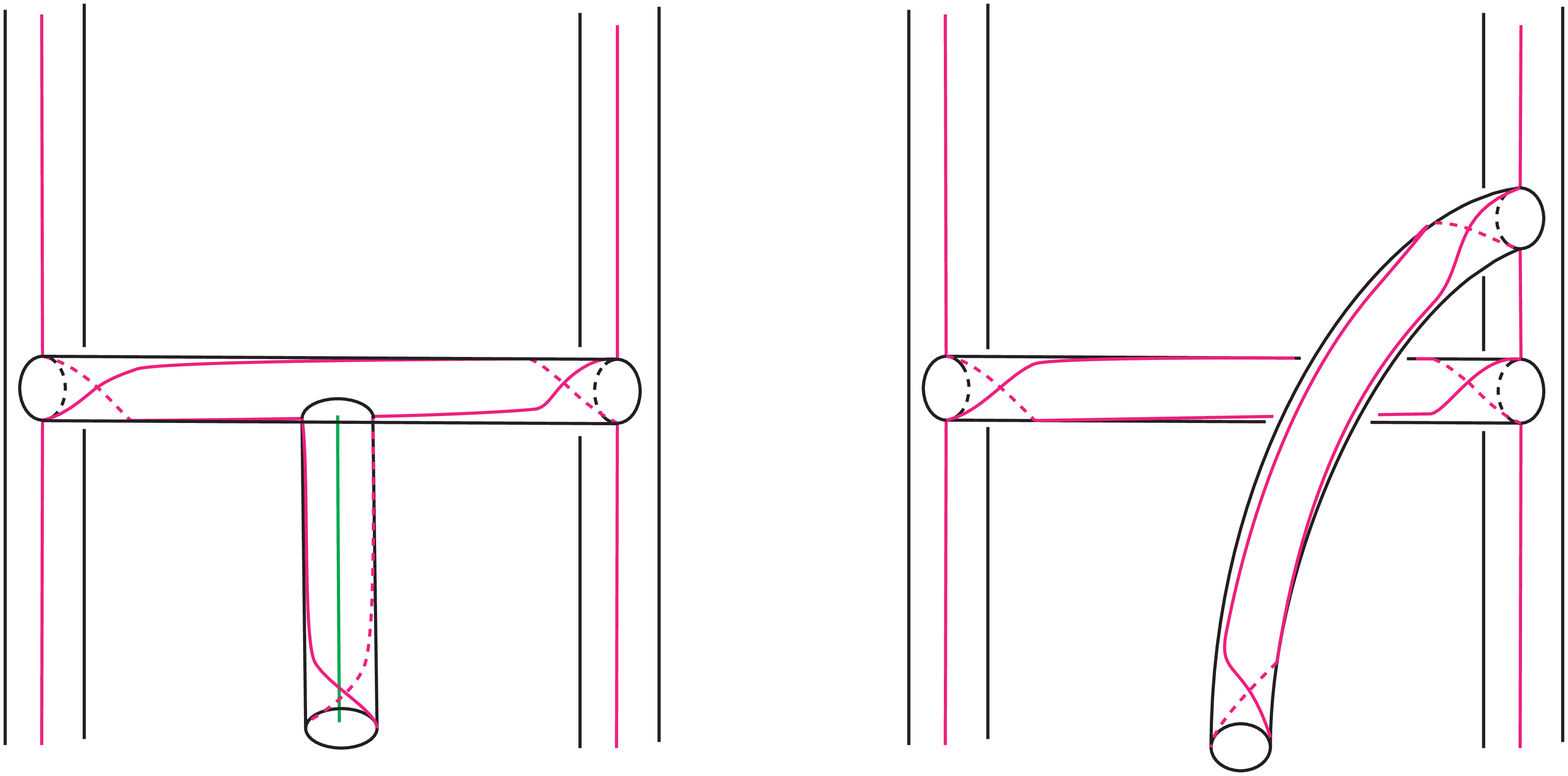}
\end{overpic}
\caption{} \label{stab}
\end{figure}

Another way to see this, without the use of contact topology, is the
following.
\begin{figure}[ht]
\begin{overpic}[width=7cm]{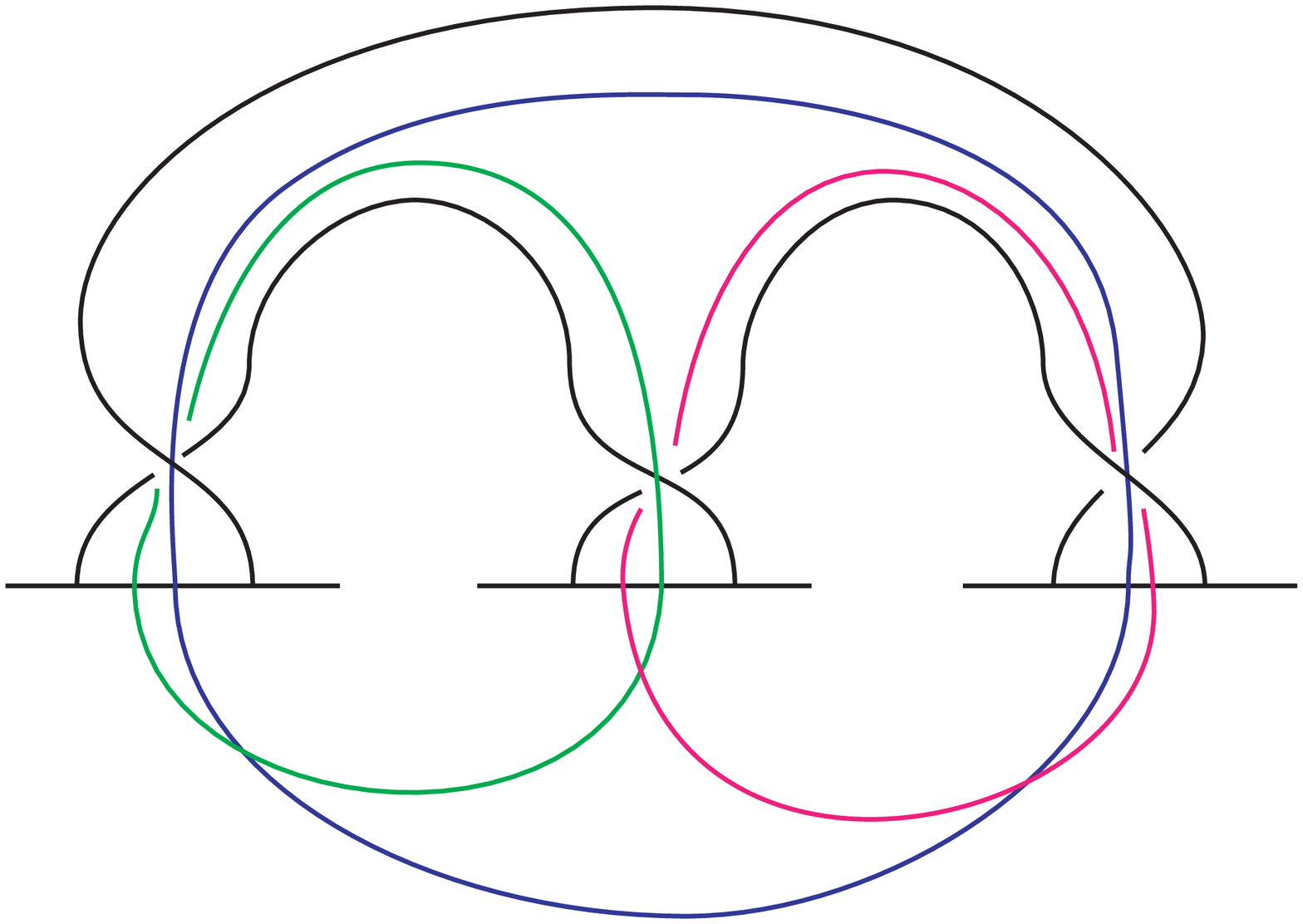}
\put(50,-3){\tiny $\gamma$} \put(30,12){\tiny
$\gamma'$}\put(68,10){\tiny $\gamma''$}
\end{overpic}
\caption{} \label{puncturedtorus}
\end{figure}
Let $c=c_1\cup c_2$, $c'=c_1\cup d$, and $c''=c_2\cup d$.  Also let
$\gamma$, $\gamma'$, and $\gamma''$ be closed curves as depicted in
Figure~\ref{puncturedtorus}, which extend $c$, $c'$, $c''$. Then
$R_{\gamma'}\circ R_{\gamma}= R_{\gamma''}\circ R_{\gamma'}$ (on a
punctured torus). (Here $R_\gamma$ is the positive Dehn twist about
$\gamma$.)
\end{proof}

\begin{prop}
Let $(S',h')$ be a stabilization of $(S,h)$ along an arc $c$ of
complexity $n$. Then $\EH(S,h)$ is mapped to $\EH(S',h')$ under the
stabilization and handleslide maps defined in \cite{OS1}.
\end{prop}

\begin{proof}
Let $(\Sigma,\beta,\alpha)$ be a compatible Heegaard splitting for
$(M,\Gamma)$ corresponding to some basis for $(S,h)$.  (Similarly
define $(\Sigma,\beta^{(i)},\alpha^{(i)})$, for $(S^{(i)},h^{(i)})$,
defined subsequently.) Suppose $(S',h')$ is a stabilization of
$(S,h)$ of complexity zero. Then the map
$\Phi:\SFH(\beta,\alpha)/\{\pm 1\}\rightarrow
\SFH(\beta',\alpha')/\{\pm 1\}$, defined as a composition of
handleslide maps and stabilization maps, sends $\EH(S,h)\mapsto
\EH(S',h')$ by Lemma~\ref{lemma: complexity zero} and the invariance
under basis change. Assume inductively that stabilizations of
complexity at most $n-1$ take $\EH$ classes to $\EH$ classes.  Then, by
Proposition~\ref{prop: stab} (we are using the same notation as the
proposition), the map $\Phi_{12}: \SFH(\beta',\alpha')/\{\pm
1\}\rightarrow \SFH(\beta'',\alpha'')/\{\pm 1\}$ sends
$\EH(S',h')\mapsto \EH(S'',h'')$ and the map $\Phi_{02}:
\SFH(\beta,\alpha)/\{\pm 1\}$ $\rightarrow$
$\SFH(\beta'',\alpha'')/\{\pm 1\}$ sends $\EH(S,h)\mapsto
\EH(S'',h'')$. Now, $\Phi_{01}:\SFH(\beta,\alpha)/\{\pm 1\}\rightarrow
\SFH(\beta',\alpha')/\{\pm 1\}$ satisfies $\Phi_{02}=\Phi_{12}\circ
\Phi_{01}$ by Theorem~2.1 of \cite{OS5}, which states that the maps
$\Phi_{ij}$ do not depend on the particular sequence of
handleslides, stabilizations, and isotopies chosen. This implies
that $\Phi_{01}$ maps $\EH(S,h)\mapsto \EH(S',h')$.
\end{proof}

\section{Properties of the contact class} \label{section: properties}

In this section we collect some basic properties of the contact
class $\EH(M,\Gamma,\xi)$.  Most of the properties are analogs of
properties of the contact class that are well-known in the case when
$M$ is closed.  The theorem which does not have an analog in the
closed case (for obvious reasons) is the restriction theorem
(Theorem~\ref{thm:restriction}).

Consider a partial open book decomposition $(S,h)$ for
$(M,\Gamma,\xi)$.  The notion of a {\em right-veering} $(S,h)$ can
be defined in the same way as in~\cite{HKM1}:  If for every $x\in
\bdry P-\Gamma$ and every properly embedded arc $a\subset P$ which
begins at $x$ and has both endpoints on $\bdry P-\Gamma$, $h(a)$ is
to the right of $a$, then we say $(S,h)$ is {\em right-veering}.

\begin{prop} \label{prop1}
If $(S,h)$ is not right-veering, then $(M,\Gamma,\xi)$ is
overtwisted. Any overtwisted $(M,\Gamma,\xi)$ admits a partial open
book decomposition $(S,h)$ which is not right-veering.
\end{prop}

\begin{proof}
The first assertion is proved in the same way as the analogous
statement in \cite{HKM1}.

For the second assertion, note that Example~1 of
Section~\ref{section: examples} gives a partial open book
decomposition of a neighborhood of an overtwisted disk with a
left-veering arc. Attach the neighborhood of a Legendrian arc $a$
which connects the boundary of the overtwisted disk to
$\Gamma\subset \bdry M$, and then complete it to a partial open book
for $(M,\Gamma)$.  The left-veering arc from Example 1 survives to
give a left-veering arc.
\end{proof}

\begin{prop} \label{prop: not rv}
If $(M,\Gamma,\xi)$ admits a partial open book decomposition $(S,h)$
which is not right-veering, then $\EH(M,\Gamma,\xi)=0$.
\end{prop}

\begin{proof}
Same as that of \cite{HKM2}.
\end{proof}

By combining Propositions~\ref{prop1} and \ref{prop: not rv}, we
obtain the following:

\begin{cor}
If $(M,\Gamma,\xi)$ is overtwisted, then $\EH(M,\Gamma,\xi)=0$.
\end{cor}

The next proposition describes the effect of Legendrian surgery on
the contact invariant.

\begin{prop}\label{prop: legsurg}
If $\EH(M,\Gamma,\xi)\not=0$ and $(M',\Gamma',\xi')$ is obtained from
$(M,\Gamma,\xi)$ by a Legendrian $(-1)$-surgery along a closed
Legendrian curve $L$, then $\EH(M',\Gamma',\xi')\not=0$.
\end{prop}

\begin{proof}
We can easily extend $L$ to a Legendrian skeleton $K$ as given in
Theorem~\ref{thm: decomposition}.  Hence there is a partial open
book decomposition $(S,R_+(\Gamma),h)$ so that $L\subset P$. Take a
basis $\{a_1,\dots,a_r\}$ so that $a_1$ intersects $L$ once and
$a_i\cap L=\emptyset$ for $i>1$. Push the $a_i$ off to obtain $b_i$
and $c_i$, as drawn in Figure~\ref{Mult}. Let
$\alpha_i=(a_i\times\{1\})\cup (a_i\times\{-1\})$,
$\beta_i=(b_i\times\{1\})\cup (h(b_i)\times\{-1\})$, and
$\gamma_i=(c_i\times\{1\})\cup (R_\gamma\circ h(c_i)\times\{-1\})$.
Here $\gamma$ is the curve $h(L)$ on the page $S$.

\begin{figure}[ht]
\begin{overpic}[width=3.8cm]{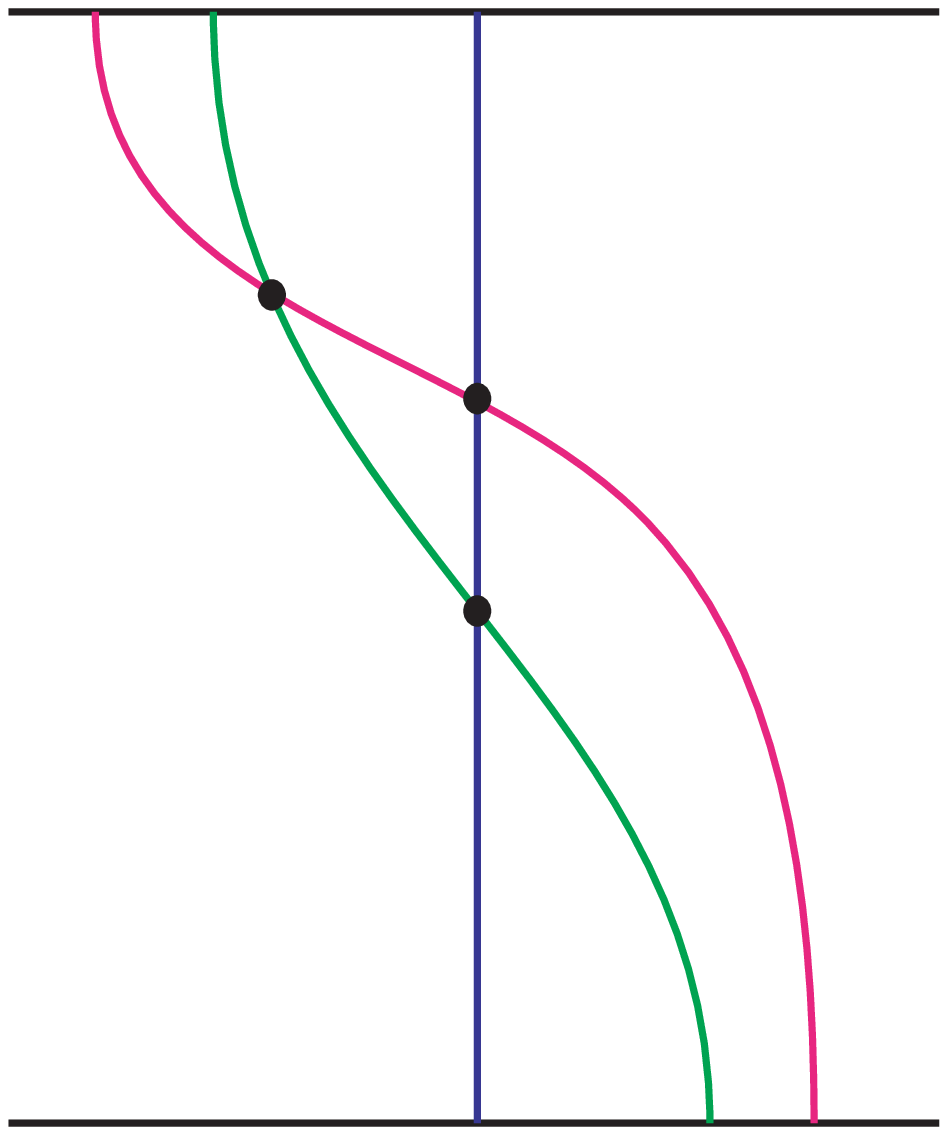}
\put(32.5,20){\tiny $a_i$} \put(52.5,15){\tiny
$b_i$}\put(73,18){\tiny $c_i$} \put(44,47){\tiny $\mathbf{x}$}
\put(15.5,66.5){\tiny $\mathbf{x}'$} \put(45,65){\tiny
$\mathbf{x}''$}
\end{overpic}
\caption{} \label{Mult}
\end{figure}

Let $\mathbf{x}\in \mathbb{T}_\beta\cap \mathbb{T}_\alpha$,
$\mathbf{x}'\in \mathbb{T}_\gamma\cap \mathbb{T}_\beta$ and
$\mathbf{x}''\in \mathbb{T}_\gamma\cap \mathbb{T}_\alpha$ be the
unique $r$-tuples which are on $S\times\{1\}$.  As Baldwin observed
in~\cite{Ba}, there is a comultiplication
$\CF(\gamma,\alpha)\rightarrow \CF(\gamma,\beta)\otimes
\CF(\beta,\alpha)$ which takes $\mathbf{x}''$ to $\mathbf{x}'\otimes
\mathbf{x}$.  We are assuming that $[\mathbf{x}]\not=0$.  Also,
$[\mathbf{x}']\not=0$ is immediate from the fact that $\beta_1$ and
$\gamma_1$ intersect only at $x'_1$, due to our choice of basis.
Therefore, we have $[\mathbf{x}'']\not=0$.
\end{proof}

In \cite{H2}, the first author exhibited a tight contact structure
on a handlebody which became overtwisted after Legendrian surgery.
By Proposition~\ref{prop: legsurg}, its contact invariant must
vanish.

\begin{thm} \label{thm:restriction}
Let $(M,\xi)$ be a closed contact 3-manifold and $N\subset M$ be a
compact submanifold $($without any closed components$)$ with convex
boundary and dividing set $\Gamma$. If $\EH(M,\xi)\not=0$, then
$\EH(N,\Gamma, \xi|_N)\not=0$.
\end{thm}

\begin{proof}
Consider the partial open book decomposition $(S, R_+(\Gamma),h)$
for $(N,\Gamma,\xi|_N)$, obtained by decomposing $N$ into $N(K)$ and
$N-N(K)$.  Now, the complement $N'=M-N$ can be similarly decomposed
into $N(K')$ and $N'-N(K')$, where $K'$ is a Legendrian graph in
$N'$ with univalent vertices on $\Gamma$.  We may assume that the
univalent vertices of $K'$ and $K$ do not intersect. Then an open
book decomposition for $(M,\xi)$ can be obtained from the Heegaard
decomposition $(N-N(K))\cup N(K')$ and $(N'-N(K'))\cup N(K)$.
Indeed, both handlebodies are easily seen to be product disk
decomposable. The page $T$ for $(M,\xi)$ can be obtained from
$(S,R_+(\Gamma),h)$ by successively attaching $1$-handles, subject
to the condition that none of the handles be attached along $\bdry
P$, where $P=S-\overline{R_+(\Gamma)}$.  The monodromy map $g:
T\rightarrow T$ extends $h: P\rightarrow S$.

Let $\{a_1,\dots, a_r, a'_1,\dots,a'_s\}$ be a basis for $T$ which
extends a basis $\{a_1,\dots,a_r\}$ for $(S,R_+(\Gamma))$. For such
an extension to exist, $T-P$ must be connected. This is possible if
$M-N$ is connected --- simply take suitable stabilizations to
connect up the components of $T-P$. If $M-N$ is disconnected, we
apply a standard contact connected sum inside $M-N$ to connect up
disjoint components of $M-N$. This has the effect of attaching
$1$-handles to $T$ away from $P$ and extending the monodromy map by
the identity. The contact manifold $(M,\xi)$ has been modified, but
it is easy to see that the contact class of the connected sum is
nonzero if and only if the original contact class is nonzero.

Let $\mathbf{x}$ be the generator of $\EH(N,\Gamma,\xi|_N)$ with
respect to $\{a_1,\dots,a_r\}$ and $(\mathbf{x},\mathbf{x'})$ be the
generator of $\EH(M,\xi)$ with respect to $\{a_1,\dots, a_r,
a'_1,\dots,a'_s\}$. If $\bdry (\sum_i c_i\mathbf{y}_i)=\mathbf{x}$,
then we claim that $\bdry (\sum_i
c_i(\mathbf{y}_i,\mathbf{x'}))=(\mathbf{x},\mathbf{x'})$. Indeed,
each of the intersection points of $\mathbf{x'}$ must map to itself
via the constant map --- this uses up all the intersection points of
$\mathbf{x'}$.  We then erase all the $\alpha_i$ and $\beta_i$
corresponding to $\mathbf{x'}$, and are left with the Heegaard
diagram for $(S,R_+(\Gamma))$.
\end{proof}

\s\n {\bf Comparison with other invariants.} We now make some
remarks on the relationship to the contact invariant in the closed
case and to Legendrian knot invariants.

\s\n 1. If we start with a closed $(M,\xi)$, then we can remove a
standard contact 3-ball $B^3$ with $\Gamma_{\bdry B^3}=S^1$. The
sutured manifold is called $M(1)$ in Juh\'asz~\cite{Ju1}.  In this
case, $\SFH(-M(1))=\HFhat(-M)$, and the contact element in
$\SFH(-M(1))$ coincides with the Ozsv\'ath-Szab\'o contact class in
$\HFhat(-M)$. (Think of the disk $R_+$ being squashed to a
point to give the basepoint $z$.)

\s\n 2. A Legendrian knot $L$ has a standard neighborhood $N(L)$
which has convex boundary. The dividing set $\Gamma_{\bdry N(L)}$
satisfies $\#\Gamma_{\bdry N(L)}=2$ and $|\Gamma_{\bdry N(L)}\cap
\bdry D^2|=2$, where $D^2$ is the meridian of $N(L)$.  Moreover, the
framing for $L$ induced by $\xi$ agrees with the framing induced by
the ribbon in $N(L)$ which contains $L$ and has boundary on
$\Gamma_{\bdry N(L)}$.

If $(M,\xi)$ is closed, then $\EH(M-N(L),\Gamma_{\bdry N(L)},\xi)$
is an invariant of the Legendrian knot $L$ which sits in
$\SFH(-(M-N(L)),-\Gamma_{\bdry N(L)})$. On the other hand, if we
choose the suture $\Gamma$ on $\bdry N(L)$ to consist of two
parallel meridian curves, then
$\SFH(-(M-N(L)),-\Gamma)=\HFKhat(-M,L)$. Hence we are slightly off
from the Legendrian knot invariants of
Ozsv\'ath-Szab\'o-Thurston~\cite{OST} and
Lisca-Ozsv\'ath-Stipsicz-Szab\'o~\cite{LOSS}
--- the Legendrian knot invariant currently does not sit in
$\HFKhat(-M,L)$. We also only have one invariant, rather than
two.

\section{Examples} \label{section: examples}

In this section we calculate a few basic examples.

\begin{figure}[ht,grid]
\begin{overpic}[width=16cm]{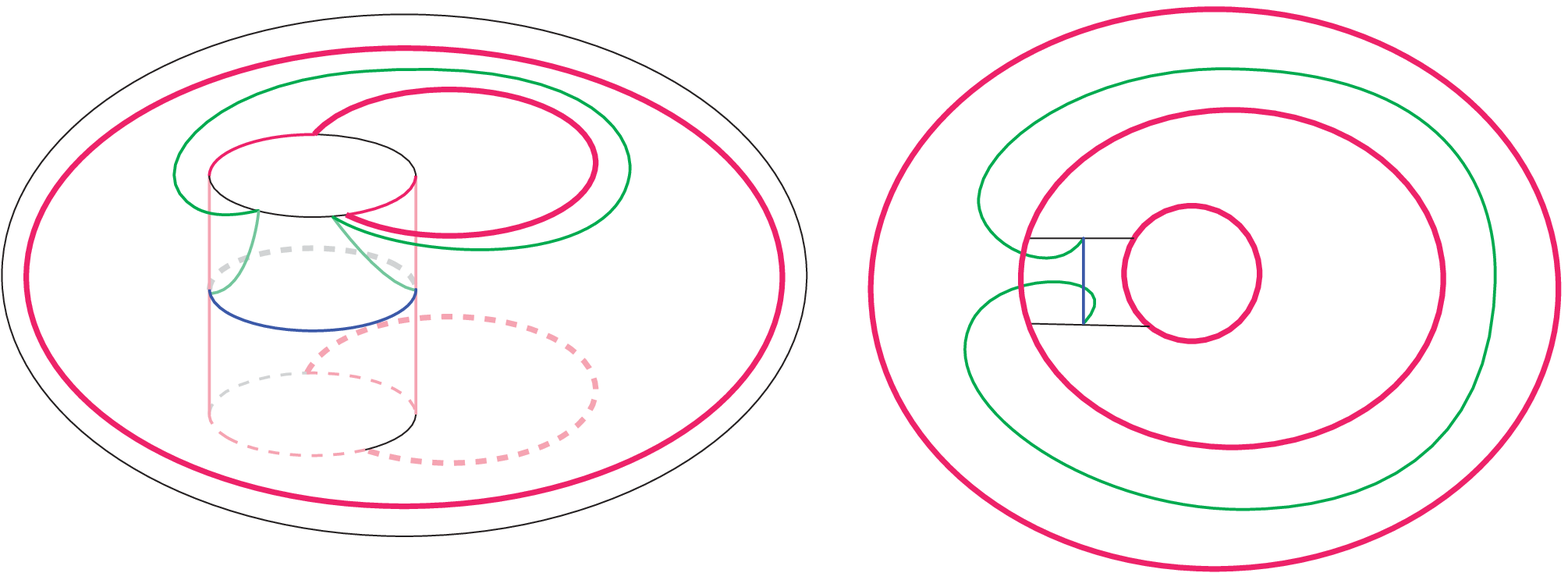}
\put(17,14.3){\tiny $a$} \put(40,23){\tiny $h(a)$} \put(43,20){\tiny
$+$} \put(30,26){\tiny $-$} \put(68.8,14.3){\tiny $x$}
\put(68.8,22){\tiny $x$} \put(69.7,18.7){\tiny $y$} \put(78,5){\tiny
$h(b)$}
\end{overpic}
\caption{} \label{OTdisk}
\end{figure}

\s\n {\bf Example 1:}  Standard neighborhood of an overtwisted disk.
Let $D$ be a convex surface which is a slight outward extension of
an overtwisted disk.  Then consider its $[0,1]$-invariant contact
neighborhood $M=D\times[0,1]$.  After rounding the edges, we obtain
$\Gamma_{\bdry M}$ which consists of three closed curves, two which
are $\gamma\times\{0,1\}$, where $\gamma$ is a closed curve in the
interior of $D^2$, and one which is $\bdry D^2\times\{{1\over 2}\}$.
Our Legendrian graph $K$ is an arc $\{p\}\times[0,1]$, where $p\in
\gamma$.  Then $M-N(K)$, after edge-rounding, is a solid torus whose
dividing set consists of two parallel homotopically nontrivial
curves, each of which intersects the meridian once. See the
left-hand diagram of Figure~\ref{OTdisk} for $M-N(K)$.  The meridian
of $M-N(K)$ gives a product disk decomposition, and the
decomposition into $M-N(K)$ and $N(K)$ gives rise to a partial open
book decomposition. The left-hand diagram of Figure~\ref{OTdisk}
shows the arc $a$ on $P$. The arc $a$ is isotoped through the
fibration $N(K)$ (rel endpoints) and then through the fibration
$M-N(K)$ (rel endpoints). The resulting arc on $S\times\{1\}$ is
$h(a)$.

Next consider the right-hand diagram of Figure~\ref{OTdisk}, which
shows a page $S$ of the partial open book decomposition.  The region
$R_+(\Gamma)\subset S$ has two connected components --- an annulus
and a disk ---  the red curves denote their boundary $\Gamma$. The
region $P$ is the $1$-handle which connects the annulus and the
disk. The blue arc is the arc $a\times\{-1\}$ and the green arc
represents $h(b)\times\{-1\}$, both viewed as sitting on
$S\times\{-1\}$. The unique intersection point $x$ between
$a\times\{1\}$ and $b\times\{1\}$ is now viewed as two points on
$S\times\{-1\}$, both labeled $x$.  (This way, all of the
holomorphic disk counting can be done on $S\times\{-1\}$.)

Now, $\CF(\beta,\alpha)$ is generated by two points $x,y$. (This is
because $h(a)$ is to the left of $a$ at one of its endpoints.) It is
easy to see that $\bdry y=x$, so $\SFH(-M,-\Gamma)=0$ and
$\EH(M,\Gamma,\xi)=0$.

\begin{figure}[ht]
\begin{overpic}[width=12cm]{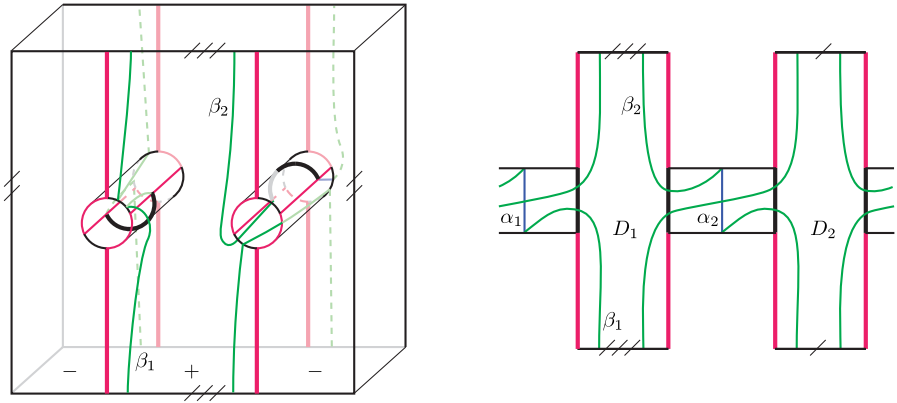}
\end{overpic} \caption{The top and bottom of the two annuli are identified.} \label{product}
\end{figure}

\s\n {\bf Example 2:}  A $[0,1]$-invariant neighborhood $M$ of a
convex surface $F$, none of whose components $F-\Gamma_F$ are disks.
Figure~\ref{product} depicts the case when $F$ is a torus with two
parallel dividing curves.  Let $\gamma_1$ and $\gamma_2$ be the two
components of $\Gamma_F$.  Then take $K=\{p_1,p_2\}\times[0,1]$,
where $p_i\in \gamma_i$, $i=1,2$.  It is not hard to see that
$M-N(K)$ is product disk decomposable.  The page $S$ is obtained
from $R_+(\Gamma)$, which is a disjoint union of two annuli, by
attaching two bands as in Figure~\ref{product}.  Hence $P$ consists
of two connected components, one with a cocore $a_1$ and the other
with a cocore $a_2$.  The black arcs on the left-hand diagram of
Figure~\ref{product} are obtained from $a_i$ by isotopy through
$N(K)$ rel endpoints, and the green arcs $h(a_i)$ by further
isotoping through $M-N(K)$. The arc $h(a_i)$ is obtained from $a_i$
via what looks like a ``half positive Dehn twist" about a boundary
component. The arc $h(a_1)$ only intersects $a_1$ and the arc
$h(a_2)$ only intersects $a_2$; hence $h(a_1),a_1$ and $h(a_2),a_2$
are independent. To determine $\SFH(-M,-\Gamma)$ and verify that
$\EH(M,\Gamma,\xi)\not=0$, the only nontrivial regions that need to
be considered are $D_1$ and $D_2$, which are annuli with corners.
(All other domains nontrivially intersect the suture $\Gamma$.) We
now apply Lipshitz' formula (Corollary~4.10 of \cite{Li1})
$$\mu(D_i)=n_{\mathbf{x}}(D_i)+n_{\mathbf{y}}(D_i)+e(D_i),$$
for computing the index of $D_i$.  Here $e(D_i)$ is the {\em Euler
measure} of $D_i$ and $n_{\mathbf{x}}(D_i)$ is the (weighted)
intersection number of the $r$-tuple $\mathbf{x}$ with $D_i$.  We
compute that $\mu(D_i)= 2({1\over 4}) +2({1\over 4}) -1=0$, which
has index $\not=1$. Hence it follows that $\SFH(-M,-\Gamma) = \Z^4=
\Z^2\otimes \Z^2$ and $\EH(M,\Gamma,\xi)\not=0$. (The fact that we
can view $\SFH(-M,-\Gamma)$ as $\Z^2\otimes \Z^2$ follows from the
observation that the two intersection points of $\alpha_1$ and
$\beta$ and the two intersection points of $\alpha_2$ and $\beta$
are independent.  We can also observe that $(M,\Gamma)$ admits an
annulus decomposition into two solid tori of the type considered in
the next example with $n=2$, and each sutured solid torus
contributes $\Z^2$. Then we can use the tensor product formula,
proved by Juh\'asz~\cite[Proposition~8.10]{Ju2}.) The general
situation is similar.

\begin{figure}[ht]
\begin{overpic}[width=12cm]{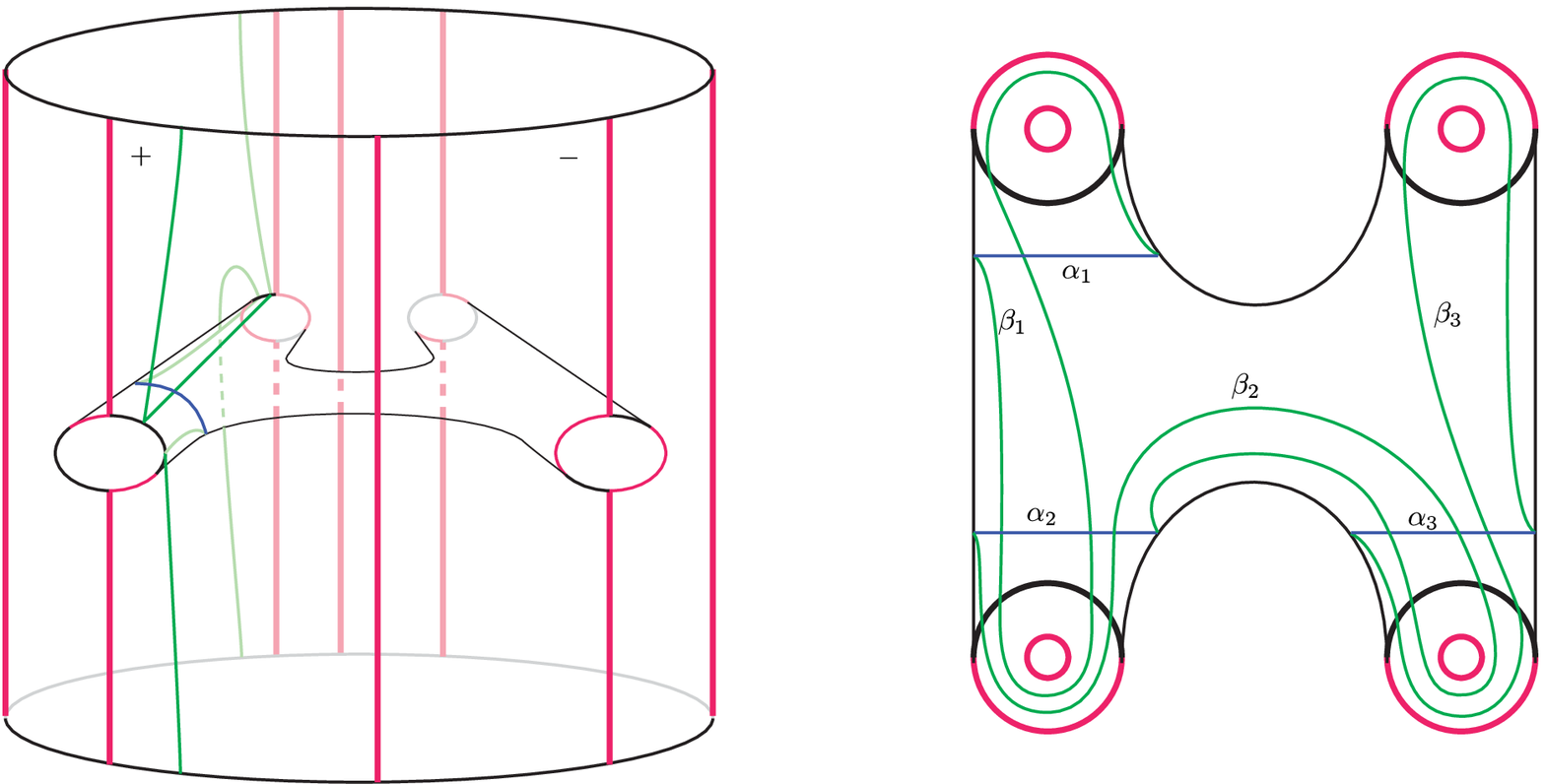}
\end{overpic} \caption{} \label{torus4}
\end{figure}

\s\n {\bf Example 3:}  Solid torus $M=S^1\times D^2$ with slope
$\infty$ and $\#\Gamma=2n>2$.  According to the classification of
tight contact structures on the solid torus \cite{Gi3,H1}, there is
a 1-1 correspondence between isotopy classes of tight contact
structures on $(M,\Gamma)$ rel boundary and isotopy classes of
dividing sets on the meridian disk $D^2=\{pt\}\times D^2$. Figure
~\ref{torus4} depicts the case $n=4$.  There are two tight contact
structures with $\bdry$-parallel $\Gamma_{D^2}$. (A convex surface
$T$ has {\em $\bdry$-parallel} dividing set $\Gamma_T$ if each
component of $\Gamma_T$ cuts off a half-disk which intersects no
other component of $\Gamma_T$.) The contact structures with
$\bdry$-parallel set are distinguished by their relative Euler
class. The monodromy map $h$ in the $\bdry$-parallel case is
calculated in Figure~\ref{torus4}. On the right-hand side of
Figure~\ref{torus4}, $a_1,a_2,a_3$ is a basis for
$(S,R_+(\Gamma),h)$, in counterclockwise order. (The $\alpha_i$ are
$\bdry (a_i\times [-1,1])$.) Label the intersections $\alpha_1\cap
\beta$ by $x_1,x_2,x_3=x_1$, $\alpha_2\cap \beta$ by
$y_1,\dots,y_5=y_1$, and $\alpha_3\cap \beta$ by
$z_1,\dots,z_5=z_1$, all in {\em clockwise} order. Starting with
$\alpha_1\cap \beta$, we find that the only valid $3$-tuples are
$(x_i,y_j,z_k)$, where $i,j,k=1,2$. There are regions of
$\Sigma-\cup_i\alpha_i-\cup_i\beta_j$ which do not intersect
$\Gamma$ --- all such regions are quadrilaterals, but use the same
$\alpha_i$ or $\beta_i$ twice, so are not valid. Hence there are no
holomorphic disks, and the boundary map is the zero map. Therefore,
$\SFH(-M,-\Gamma)=\Z^8=\Z^2\otimes\Z^2\otimes \Z^2$ and
$\EH(M,\Gamma,\xi)\not=0$. Moreover, if we split $\SFH(-M,-\Gamma)$
according to relative Spin$^c$-structures, then we have a direct sum
$\Z\oplus \Z^3\oplus \Z^3\oplus \Z$. The two $\bdry$-parallel tight
contact structures use up the first and last $\Z$ summands, and the
others have $\EH$ classes that live in the remaining $\Z^3\oplus
\Z^3$. We leave it as an exercise to verify that the $\EH$ classes
of the remaining tight contact structures are nonzero. There are
$4+2$ tight contact structures $\xi_i$, $i=1,\dots,6$, corresponding
to each $\Z^3$ summand, but we expect the $\EH$ classes to
distinguish them --- this means that we believe the $\EH(\xi_i)$ are
linearly dependent.

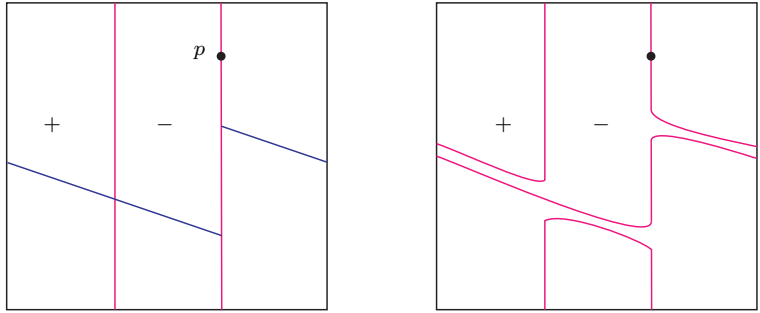
\begin{figure}[ht]
\begin{overpic}[width=10cm]{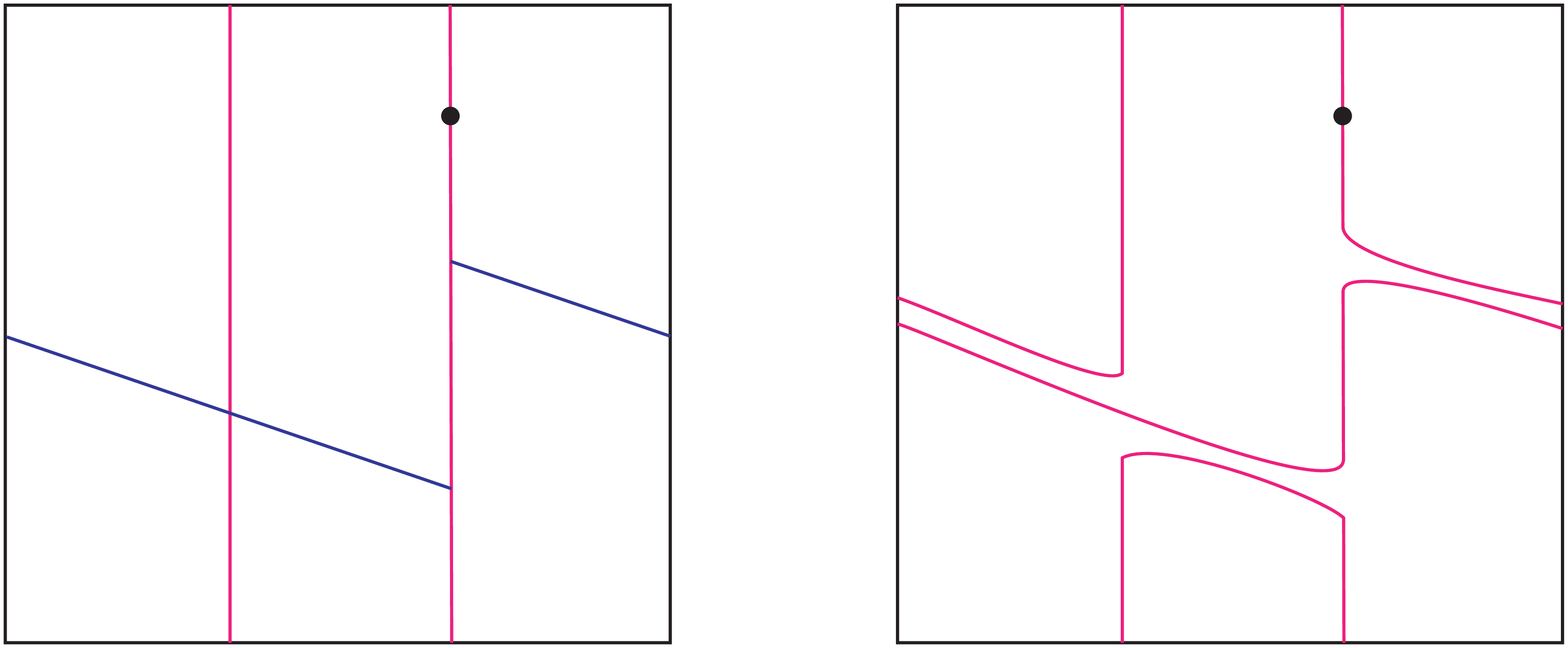}
\put(25,34){\tiny $p$} \put(20,24){\tiny $-$} \put(5,24){\tiny $+$}
\put(78,24){\tiny $-$} \put(65,24){\tiny $+$}
\end{overpic} \caption{The left-hand
diagram gives $T_0$ and the arc of attachment for the bypass, which
is attached from the front.   The right-hand diagram gives $T_1$.
Here the sides are identified and the top and the bottom are
identified to give $T^2$. The other basic slice is obtained by
switching the $+$ and $-$ signs.} \label{basicslice2}
\end{figure}

\s\n {\bf Example 4:} A basic slice $M=T^2\times[0,1]$
(see~\cite{H1}). Write $T_i=T^2\times\{i\}$, $i=0,1$, and take an
oriented identification $T^2\simeq \R^2/\Z^2$.  Normalize so that
$\#\Gamma_{T_i}=2$, $\Gamma_{T_i}$ are linear, and
$\mbox{slope}(\Gamma_{T_1})=0$, $\mbox{slope}(\Gamma_{T_0})=\infty$.
A basic slice can be obtained from the convex torus $T_0$ by
attaching a single bypass along a linear arc of slope $-1<s<0$ and
thickening. Pick a point $p$ on the connected component of
$\Gamma_{T_0}$ which contains the endpoints of the bypass arc of
attachment, as indicated in Figure~\ref{basicslice2}. Then we let
$K$ be $p$ times an interval. It is not hard to see that $M-N(K)$ is
product disk decomposable.

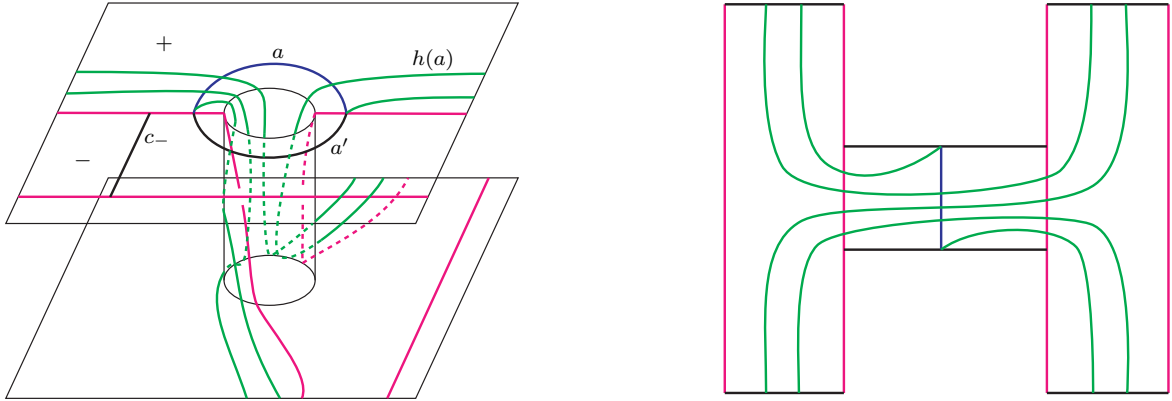
\begin{figure}[ht]
\begin{overpic}[width=15.5cm]{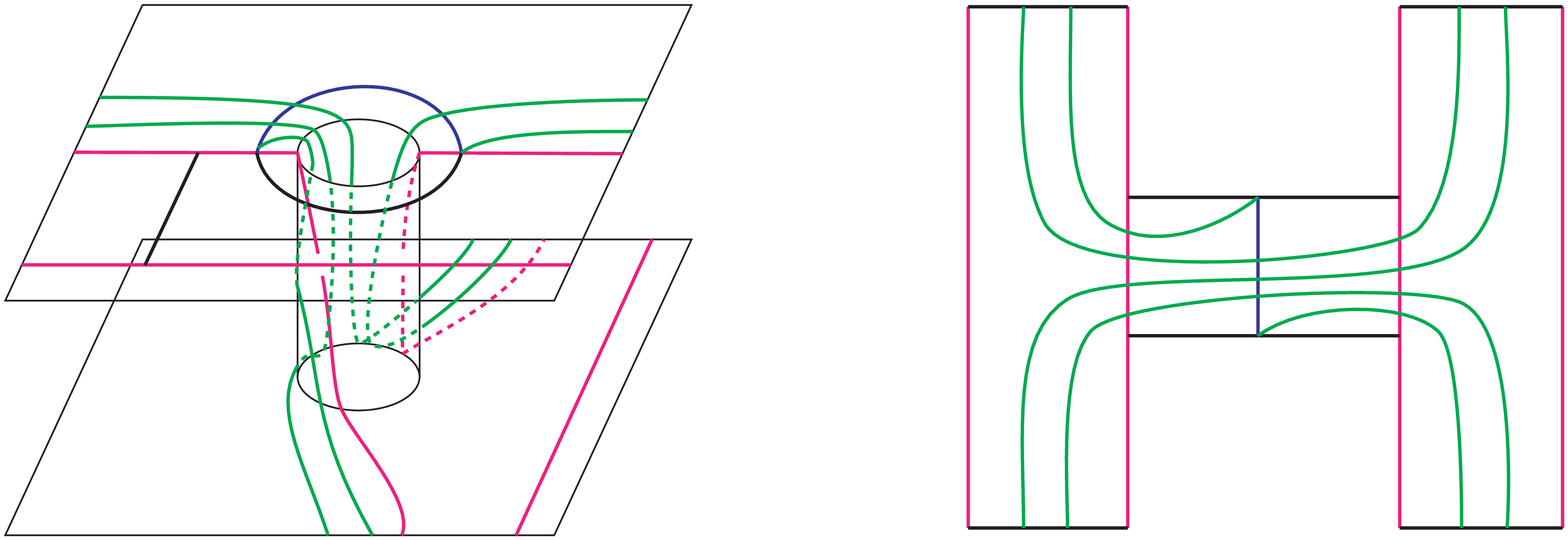}
\put(13,30){\tiny $+$} \put(6,20){\tiny $-$} \put(23,29.3){\tiny
$a$} \put(35,28.65){\tiny $h(a)$} \put(28,21){\tiny $a'$}
\put(12,22){\tiny $c_-$}
\end{overpic} \caption{The left-hand diagram calculates the monodromy map $h$.
The arc $a$ is pushed across $N(K)$ to give $a'$.  Then $a'$ is
pushed across $M-N(K)$ to give $h(a)$. The right-hand diagram gives
the monodromy map for a basic slice. Here the top and the bottom of
the two annuli are identified.} \label{basicslice}
\end{figure}

The page $S$ is a thrice-punctured sphere obtained from
$R_+(\Gamma)$, a disjoint union of two annuli, by adding a 1-handle
to connect the two annuli. If $a\subset P$ is the cocore of the
1-handle, then we compute that $h$ is a positive Dehn twist about
the connected component of $\bdry S$ which contains the endpoints of
$a$.  See Figure~\ref{basicslice}.   One easily computes that
$\SFH(-M,-\Gamma)=\Z^4$ and $\EH(M,\Gamma,\xi)\not=0$. There are two
basic slices, which are distinguished by the relative Euler class.
They account for $\Z\oplus \Z\subset \SFH(-M,-\Gamma)$. The
remaining $\Z\oplus \Z$ come from tight contact structures with an
extra ${\pi}$-rotation, and are also distinguished by their relative
Euler classes. (See Example 6.)

\begin{figure}[ht]
\begin{overpic}[width=8cm]{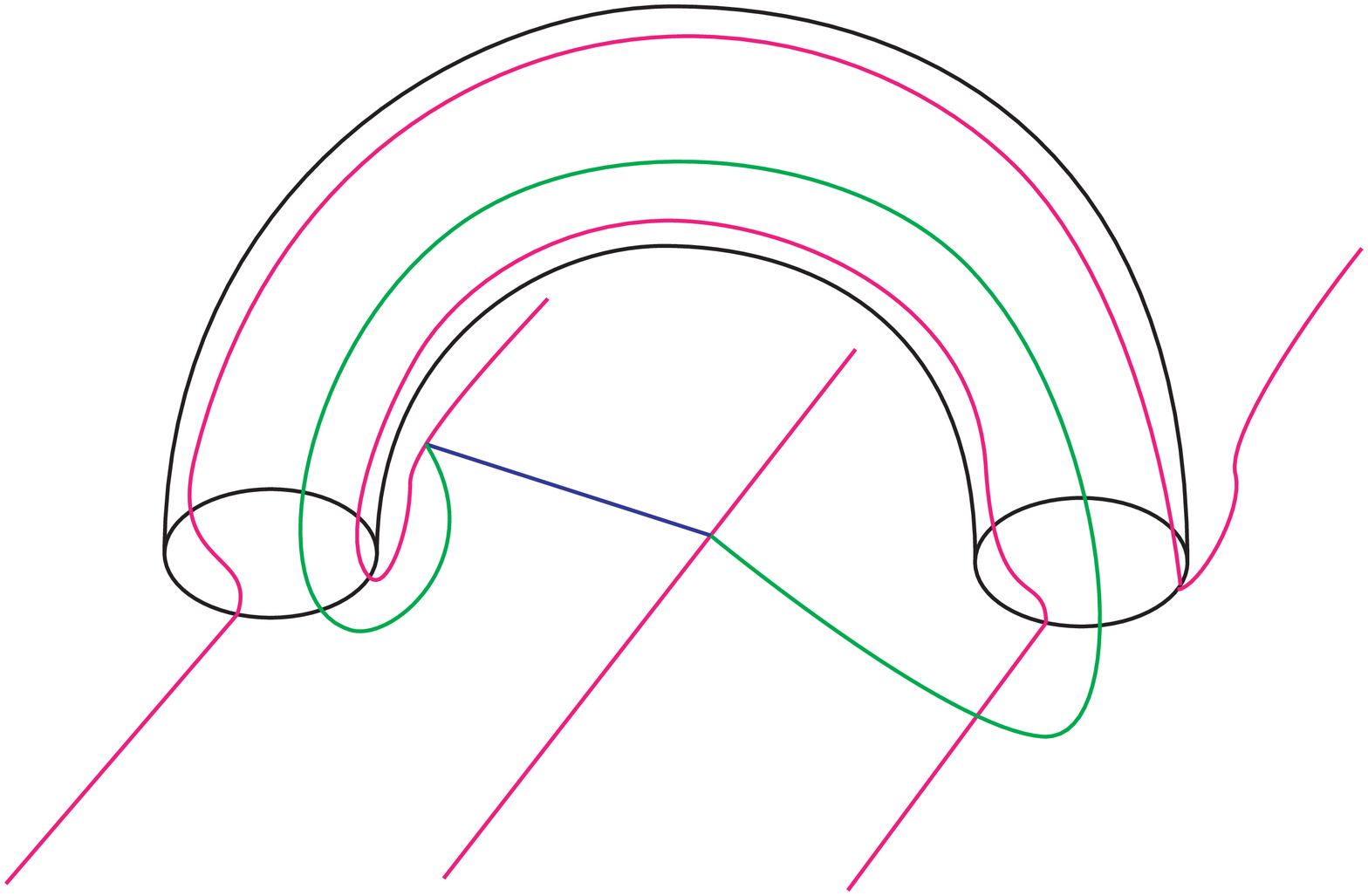}
\put(41.5,30.5){\tiny $a$} \put(25,16){\tiny $a'$} \put(56,17){\tiny
$a''$}\put(22,10){\tiny $+$} \put(53,10){\tiny $-$}
\end{overpic} \caption{} \label{bypass}
\end{figure}

\begin{figure}[ht]
\begin{overpic}[width=11.5cm]{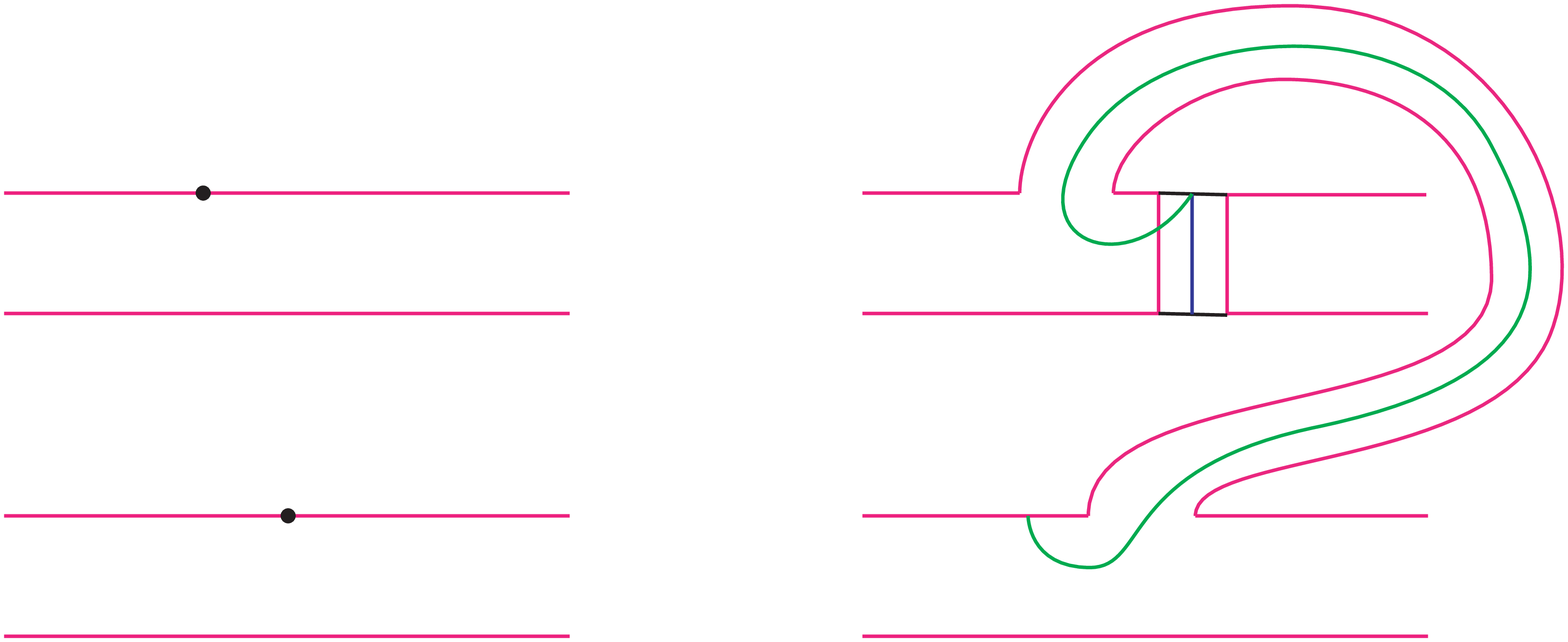}
\put(18,24){\tiny $+$} \put(18,3){\tiny $+$} \put(63,24){\tiny $+$}
\put(63,3){\tiny $+$} \put(12,30){\tiny $p_1$} \put(17.5,9.8){\tiny
$p_2$}
\end{overpic} \caption{The left-hand diagram shows a portion of $S$ before bypass attachment and
the right-hand diagram is $S'$, after modification.  The green arc
is $a'$ and the blue arc is $a$.} \label{effectofbypass}
\end{figure}
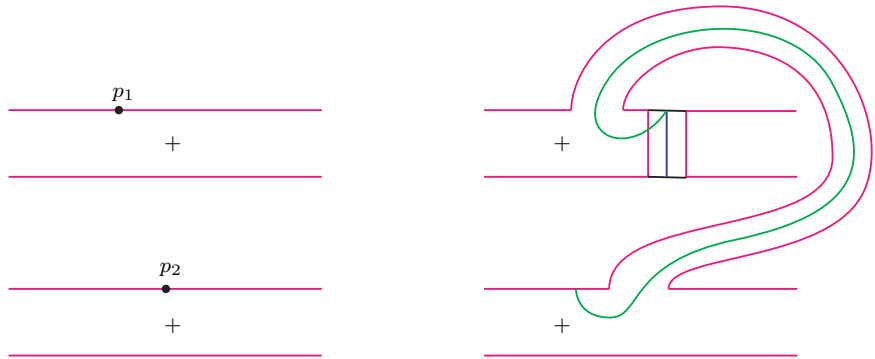

\s\n {\bf Example 5:} Attaching a bypass.  Consider the contact
structure $(M,\Gamma,\xi)$.  Let $D$ be a bypass which is attached
to $M$ from the exterior, along a Legendrian arc of attachment
$c\subset \bdry M$.  If $M'=M\cup N(D)$ and $(M',\Gamma',\xi')$ is
the resulting contact 3-manifold with convex boundary, then we
express the monodromy map $(S',R_+(\Gamma'),h')$ for
$(M',\Gamma',\xi')$ in terms of $(S,R_+(\Gamma),h)$ for
$(M,\Gamma,\xi)$.  Assume without loss of generality that $c$ does
not intersect the neighborhood $N(K)$ of the Legendrian skeleton $K$
for $M$.  Let $p_1, p_2$ be the endpoints of $c$ on $\Gamma$ and let
$c_\pm=c\cap \overline{R_\pm(\Gamma)}$. Next let $\bdry D=c\cup d$,
where $c$ and $d$ intersect only at their endpoints $p_1,p_2$. The
key observation is that a bypass attachment corresponds to two
handle attachments: a 1-handle $N(d)$ attached at $p_1$ and $p_2$,
followed by a canceling 2-handle.  Then $(M-N(K))\cup N(d)$ is
product disk decomposable and gives a fibration structure
$S'\times[-1,1]$.  Let $a$ be the arc in $P'$ given in
Figure~\ref{bypass}.  Then, after isotoping it through the 2-handle,
$a$ is isotopic rel endpoints to $a'\cup a''$, where $a'$ lies in
the positive region and $a''$ lies in the negative region of $\bdry
((M-N(K))\cup N(d))$.  Next, isotop $a''$ through the fibration
$M-N(K)$ to obtain $a'''$ which lies in the positive region of
$\bdry (M-N(K))$.  Then $h'(a)=a'\cup a'''$. All the other arcs of
$P'$, i.e., those that were in $P$, are unaffected. Summarizing,
$S'$ is obtained from $S$ by attaching a 1-handle from $p_1$ to
$p_2$ as in Figure~\ref{effectofbypass}, the new $P'$ is the union
of $P$ and a neighborhood of the arc $a$, $h'(a)=a'\cup a'''$, and
$h'=h$ on $P$.

\s\n {\bf Example 6:} Tight contact structures on $T^2\times I$.
Using Examples 4 and 5, we analyze partial open book decompositions
obtained from a basic slice by successively attaching bypasses.
These calculations were motivated by computations of open books for
tight contact structures on $T^3$ by Van Horn-Morris~\cite{VH}, and
echo the results obtained there. In particular, we will build up to
an open book for the following tight contact structure on $T^2\times
[0,4]$:  If $(x,y,t)$ are coordinates on $T^2\times [0,4]$, then
consider $\ker \alpha$, where $\alpha= \cos ({\pi\over 2}t) dx -\sin
({\pi\over 2}t) dy$.  Perturb $\ker \alpha$ along $\bdry
(T^2\times[0,4])$ so that so it becomes convex with
$\#\Gamma_{T_0}=\#\Gamma_{T_4}=2$, and
$\mbox{slope}(\Gamma_{T_0})=\mbox{slope}(\Gamma_{T_4})=\infty$.

\begin{figure}[ht]
\begin{overpic}[width=14cm]{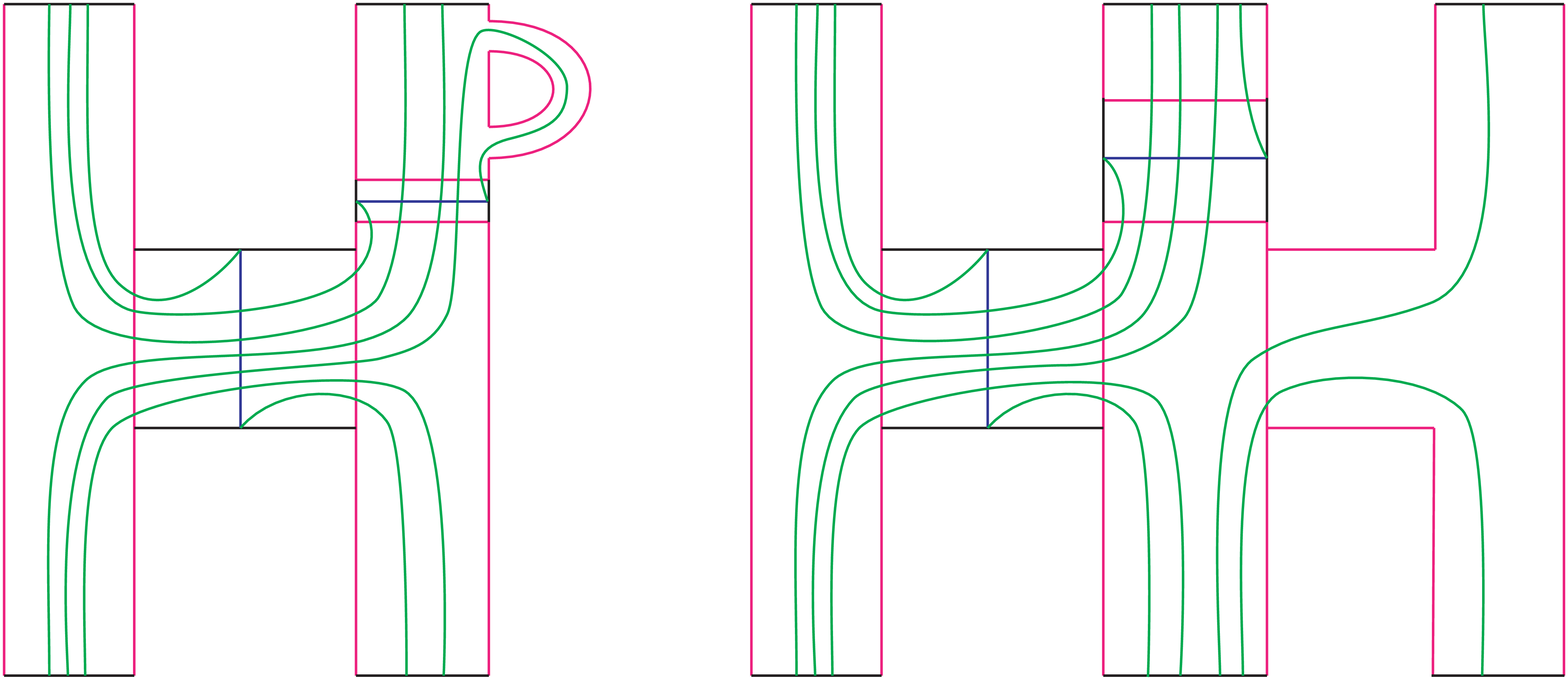}
\end{overpic} \caption{The monodromy map.  The right-hand diagram is
the same as the left-hand one, but drawn more symmetrically. The top
and bottom of each annulus is identified.} \label{tight1}
\end{figure}

\begin{figure}[ht]
\begin{overpic}[width=10cm]{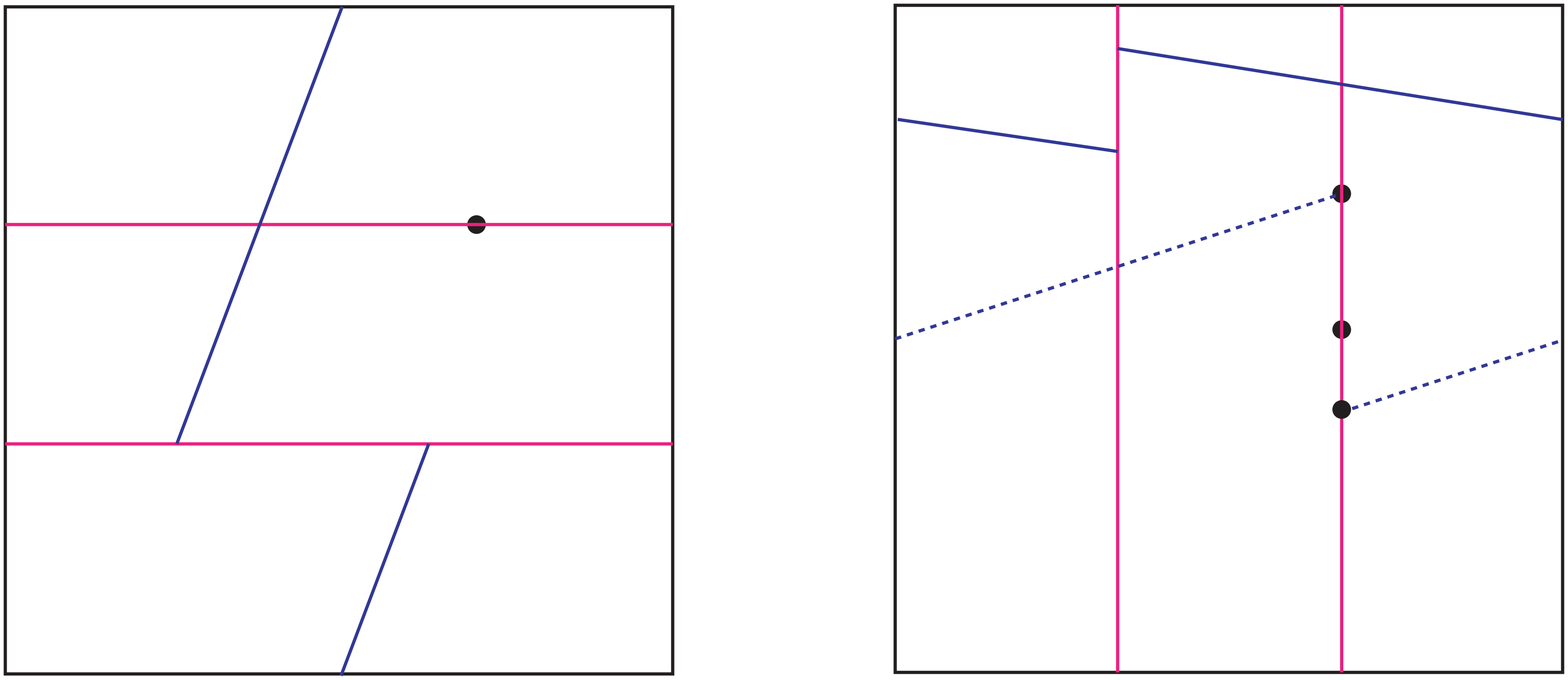}
\end{overpic} \caption{The left-hand diagram is $T_1$ with a bypass arc from the front.
After attaching the bypass we get $T_2$ on the right-hand side, with
the corresponding ``anti-bypass'' to the back, in dashed lines. The
solid lines give the bypass attached from the front to give $T_3$.}
\label{tight2}
\end{figure}

\s\n (a) Start with the basic slice $M=T^2\times[0,1]$ from Example
4, depicted in Figure~\ref{basicslice2}. Take a linear Legendrian
arc $c$ of slope $1<s<\infty$ so that the component of
$\Gamma_{T_1}$ that contains the endpoints of $c$ is different from
that containing $p$. (If the endpoints of $c$ are on the same
component as $p$, then the resulting contact structure is
overtwisted.) See the left-hand diagram of Figure~\ref{tight2}.
Attaching a bypass along $c$ and thickening yields
$M'=T^2\times[0,2]$ with $\# \Gamma_{T_i}=2$,
$\mbox{slope}(\Gamma_{T_i})=0$, $i=0,2$. Moreover, the contact
structure is obtained from the restriction of $\alpha$ to
$T^2\times[0,2]$, by perturbing the boundary.  Hence it has
``torsion'' $\approx {\pi}$, where we indicate that the term is used
informally by the use of quotations.  If $c=c_+\cup c_-$, where
$c_\pm=c\cap \overline{R_\pm(\Gamma_{T_1})}$, then we obtain $a$ (as
in Example 5) by pushing $c_+$ slightly to the right. We isotop
$c_-$ to $c_+'$ in the positive region of $\bdry (M -N(p))$
--- this is readily done by referring to Figure~\ref{basicslice}.
The result is $h'(a)=c_+c_+'$ as shown in Figure~\ref{tight1}.
The monodromy $h'$ is the restriction of $R_{\gamma_1} R_{\gamma_2}
R_{\delta_1}^{-2}$ to $P$, where $\gamma_1$ and $\gamma_2$ are
parallel to the two boundary components in the middle and $\delta_1$
is the core curve of the middle vertical annulus. (Refer to the
right-hand diagram.) Since the only regions of
$\Sigma-\cup_i\alpha_i-\cup_i\beta_i$ which do not intersect the
suture $\Gamma'$ for $M'$ are quadrilaterals, it is easy to compute
that $\EH\not=0$.  The contact class lives in
$\SFH(-M',-\Gamma')=\Z^2\otimes \Z^2$ from Example 2 (but is in a
different $\Z$ summand).

\s\n (b) Next attach another bypass along $T_2$ to obtain
$M''=T^2\times[0,3]$, where $\#\Gamma_{T_0}=\#\Gamma_{T_3}=2$,
$\mbox{slope}(\Gamma_{T_0})=0$, and $\mbox{slope}(\Gamma_{T_3})=
\infty$. The contact structure has ``torsion" $\approx {3\pi\over
2}$. See Figure~\ref{tight2}.
\begin{figure}[ht]
\begin{overpic}[width=15cm]{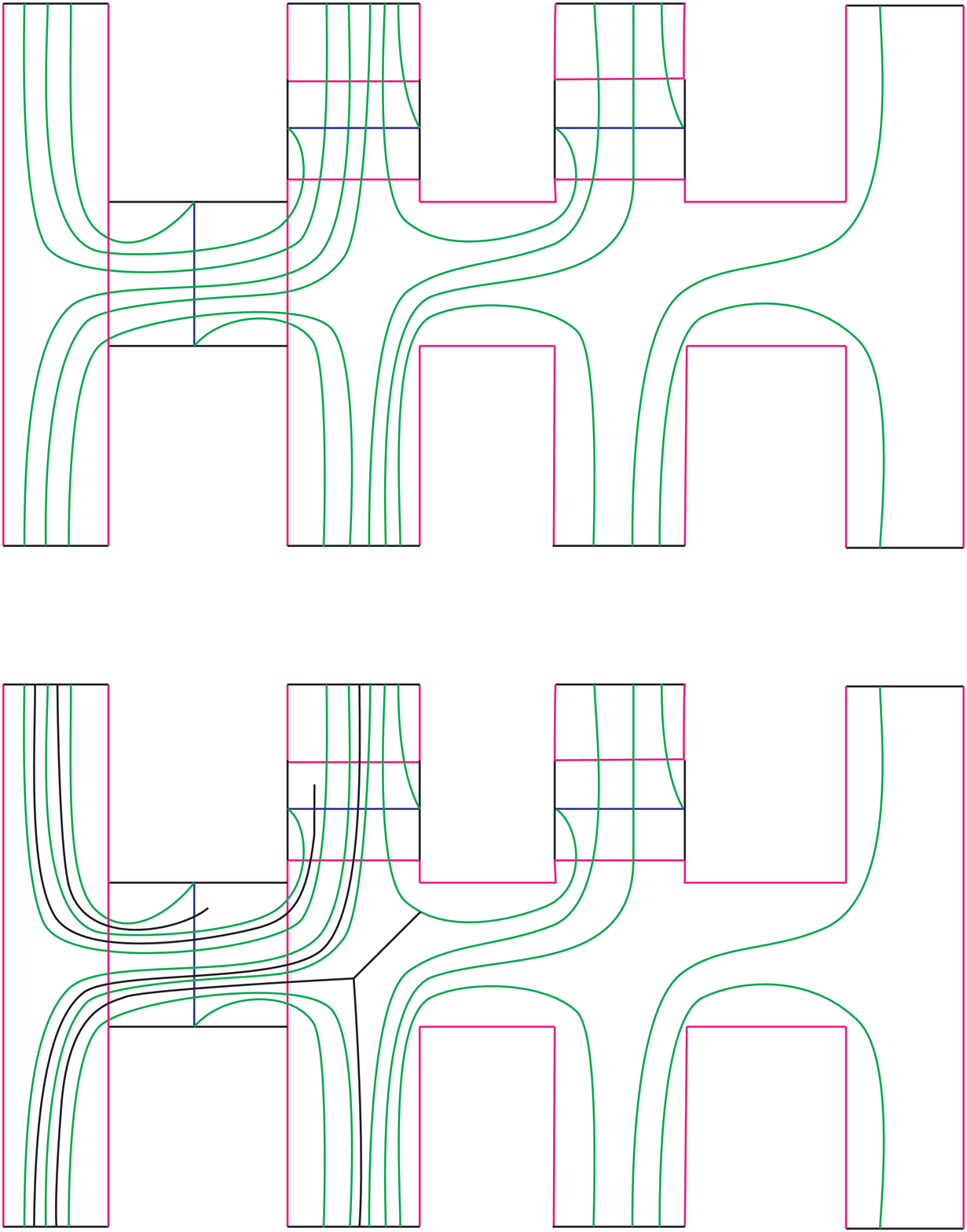}
\put(16,81.3){\tiny $a_1$} \put(24,90.6){\tiny $a_2$}
\put(46,90.6){\tiny $a_3$} \put(29.7,23){\tiny $G$}
\put(15,29){\tiny $x_1$} \put(13,15){\tiny $x_{15}=x_1$}
\put(21.3,34){\tiny $y_1$} \put(34.5,34){\tiny $y_{10}=y_1$}
\put(43.3,34){\tiny $z_1$} \put(56,34){\tiny $z_4=z_1$}
\end{overpic} \caption{}
\label{tight3}
\end{figure}
Doing a similar calculation as before, we obtain
Figure~\ref{tight3}. The monodromy $h''$ is the product
$R_{\gamma_1} R_{\gamma_2} R_{\gamma_3} R_{\delta_1}^{-2}
R_{\delta_2}^{-2}$, where $\gamma_1, \gamma_2, \gamma_3$ are
parallel to the two boundary components in the middle and $\delta_1,
\delta_2$ are the core curves of the middle vertical annuli. Let
$\{a_1, a_2, a_3\}$ be the basis, ordered from left to right in
Figure~\ref{tight3}. (The arcs are on $S''\times\{-1\}$, where $S''$
is the page. The $\alpha$-curves are
$\alpha_i=\bdry(a_i\times[-1,1])$.) To determine whether $\EH\not=0$,
we apply the Sarkar-Wang algorithm~\cite{SW} (also see \cite{Pl}) to
isotop $\beta_i$ so that the only regions of
$\Sigma-\cup_i\alpha_i-\cup_i\beta_i$ which do not intersect the
suture are quadrilaterals. (Such a diagram we call {\em
combinatorial}.)  Let $G$ (in black) be the trivalent graph in the
bottom diagram of Figure~\ref{tight3}. By isotoping $\beta_3$ across
$G$, the Heegaard diagram becomes combinatorial. Now label the
intersections $\alpha_1\cap \beta$ as $x_1,\dots,x_{14},x_{15}=x_1$
from top to bottom, $\alpha_2\cap \beta$ as
$y_1,\dots,y_9,y_{10}=y_1$ from left to right, and $\alpha_3\cap
\beta$ as $z_1,\dots,z_4=z_1$ from left to right. Observe that each
intersection of the graph with $\alpha_i$ corresponds to two
consecutive intersections with the new $\beta_3$.

We now directly prove that $\EH\not=0$.  If $c_1c_2c_3c_4$ is a
quadrilateral which avoids the suture, and $\bdry
(c_1,c_3,\mathbf{y})= (c_2,c_4,\mathbf{y})+\dots$ (for some
$\mathbf{y}$), then we say that $c_1,c_3$ are the ``from'' corners
of the quadrilateral and $c_2,c_4$ are the ``to'' corners of the
quadrilateral.  One readily calculates that the only quadrilateral
where the ``to'' corners are in $\{x_1,y_1,z_1\}$ is $y_1y_2z_1z_2$.
Now,
$$\bdry (x_1,y_2,z_2)= (x_1,y_1,z_1)+ (x_1,y_3,z_2),$$
where the second term comes from a bigon. Next, the only
quadrilateral where the ``to'' corners are in $\{x_1,y_3,z_2\}$ is
$x_1x_2y_3y_4$. We compute that
$$\bdry (x_2,y_4,z_2)= (x_1,y_3,z_2)+ (x_3,y_4,z_2).$$  Continuing,
the only quadrilateral with ``to'' corners in $\{x_3,y_4,z_2\}$ is
$x_3x_4z_3z_2$.  We have $$\bdry (x_4,y_4,z_3)=
(x_3,y_4,z_2)+(x_7,y_1,z_3).$$ There is nothing else with ``to''
corners in $\{x_7,y_1,z_3\}$.  Hence $\EH\not=0$.

\s\n (c) Now consider $M'''=T^2\times[0,4]$, where
$\#\Gamma_{T_0}=\#\Gamma_{T_4}=2$,
$\mbox{slope}(\Gamma_{T_0})=\mbox{slope}(\Gamma_{T_4})=0$, and the
contact structure $\xi'''$ has ``torsion'' $\approx 2\pi$. (This is
the contact structure described in the first paragraph of
Example~6.) The page $S'''$ is obtained from the page $S''$ for
$M''$ by attaching an annulus along the right-hand boundary in
Figure~\ref{tight3}. (We can draw the annulus to be vertical and
parallel to the previous annuli which were successively attached,
starting with $S$.)  Then the monodromy $h'''=R_{\gamma_1}
R_{\gamma_2} R_{\gamma_3} R_{\gamma_4} R_{\delta_1}^{-2}
R_{\delta_2}^{-2}R_{\delta_3}^{-2}$, where $\gamma_1, \gamma_2,
\gamma_3,\gamma_4$ are parallel to the three boundary components in
the middle and $\delta_1,\delta_2,\delta_3$ are the core curves of
the middle vertical annuli. One can probably directly show that
$\EH(M''',\Gamma''',\xi''')\not=0$ with some patience, but instead
we invoke Theorem~\ref{thm:restriction} and observe that
$(M''',\xi''')$ can be embedded in the unique Stein fillable contact
structure on $T^3$, which has nonvanishing contact invariant.

\s\n (d) Recall that a contact manifold $(M,\xi)$ has {\em $n\pi
$-torsion} with $n$ a positive integer if it admits an embedding
$(T^2\times[0,1],\eta_{n\pi })\hookrightarrow (M,\xi)$, where
$(x,y,t)$ are coordinates on $T^2\times [0,1]\simeq
\R^2/\Z^2\times[0,1]$ and $\eta_{n\pi}=\ker (\cos (n\pi t) dx - \sin
(n\pi t) dy)$. In a subsequent paper~\cite{GHV}, Ghiggini, the first
author, and Van Horn-Morris prove that if $(M,\xi)$ has torsion
$\geq 2\pi$ with $M$ closed, then $\EH(M,\xi)=0$ when
$\Z$-coefficients are used.

\section{Relationship with sutured manifold decompositions}
\label{section: sutured}

A sutured manifold decomposition $(M,\Gamma)\stackrel
{T}\rightsquigarrow (M',\Gamma')$ is {\em well-groomed} if for every
component $R$ of $\bdry M- \Gamma$, $T\cap R$ is a union of parallel
oriented nonseparating simple closed curves  if $R$ is nonplanar and
arcs if $R$ is planar.  According to \cite{HKM3}, to each
well-groomed sutured manifold decomposition there is a corresponding
convex decomposition $(M,\Gamma)\stackrel {(T,\Gamma_T)}
\rightsquigarrow (M',\Gamma')$, where $T$ is a convex surface with
Legendrian boundary and $\Gamma_T$ is a dividing set which is
$\bdry$-parallel, i.e., each component of $\Gamma_T$ cuts off a
half-disk which intersects no other component of $\Gamma_T$.

In \cite{Ju2}, Juh\'asz proved the following theorem:

\begin{thm}[Juh\'asz]
Let $(M,\Gamma)\stackrel {T}\rightsquigarrow (M',\Gamma')$ be a
well-groomed sutured manifold decomposition.  Then
$\SFH(-M',-\Gamma')$ is a direct summand of $\SFH(-M,-\Gamma)$.
\end{thm}
More precisely, $\SFH(-M',-\Gamma')$ is the direct sum
$\oplus_{\mathfrak{s}} \SFH(-M,-\Gamma,\mathfrak{s})$, where
$\mathfrak{s}$ ranges over all relative Spin$^c$-structures which
evaluate maximally on $T$. (The are called ``outer'' in \cite{Ju2}.)

In this section we give an alternate proof from the
contact-topological perspective.  We remark that our proof does not
require the use of the Sarkar-Wang algorithm. In particular, our
proof will imply the following:

\begin{thm} \label{thm: gluing}
Let $(M,\Gamma,\xi)$ be the contact structure obtained from
$(M',\Gamma',\xi')$ by gluing along a $\bdry$-parallel
$(T,\Gamma_T)$. Under the inclusion of $\SFH(-M',-\Gamma')$ in
$\SFH(-M,-\Gamma)$ as a direct summand, $\EH(M',\Gamma',\xi')$ is
mapped to $\EH(M,\Gamma,\xi)$.
\end{thm}

\begin{cor}\label{cor}
$\EH(M',\Gamma',\xi')\not=0$ if and only if $\EH(M,\Gamma,\xi)\not=0$.
\end{cor}

\begin{proof}
The ``if'' direction is given by Theorem~\ref{thm:restriction}. The
``only if'' direction follows from Theorem~\ref{thm: gluing}.
\end{proof}

The corollary is a gluing theorem for tight contact structures which
are glued along a $\bdry$-parallel dividing set, and does not
require any universally tight condition which was needed for its
predecessors, e.g., Colin's gluing theorem \cite{Co1}.

\begin{proof}
Let $T$ be a convex surface with Legendrian boundary and $\Gamma_T$
be its dividing set, which we assume is $\bdry$-parallel. Let
$(M,\Gamma,\xi)$ be a contact structure obtained from
$(M',\Gamma',\xi')$ by gluing along $(T,\Gamma_T)$. Suppose without
loss of generality that $T$ is oriented so that the disks cut off by
the $\bdry$-parallel arcs are negative regions of $T-\Gamma_T$.  Let
$L\subset T$ be a Legendrian skeleton of the positive region, with
endpoints on $\Gamma$.  We observe that there exist compressing
disks $D_i$ in $M-N(L)$ with $\#(\Gamma_{\bdry (M-N(L))}\cap \bdry
D_i)=2$ so that splitting $M-N(L)$ along the $D_i$ gives
$(M',\Gamma')$.  (The disks are basically the disks in $T$ cut off
by $L$.) See the left-hand side of Figure~\ref{mult}.
\begin{figure}[ht]
\begin{overpic}[width=11.5cm]{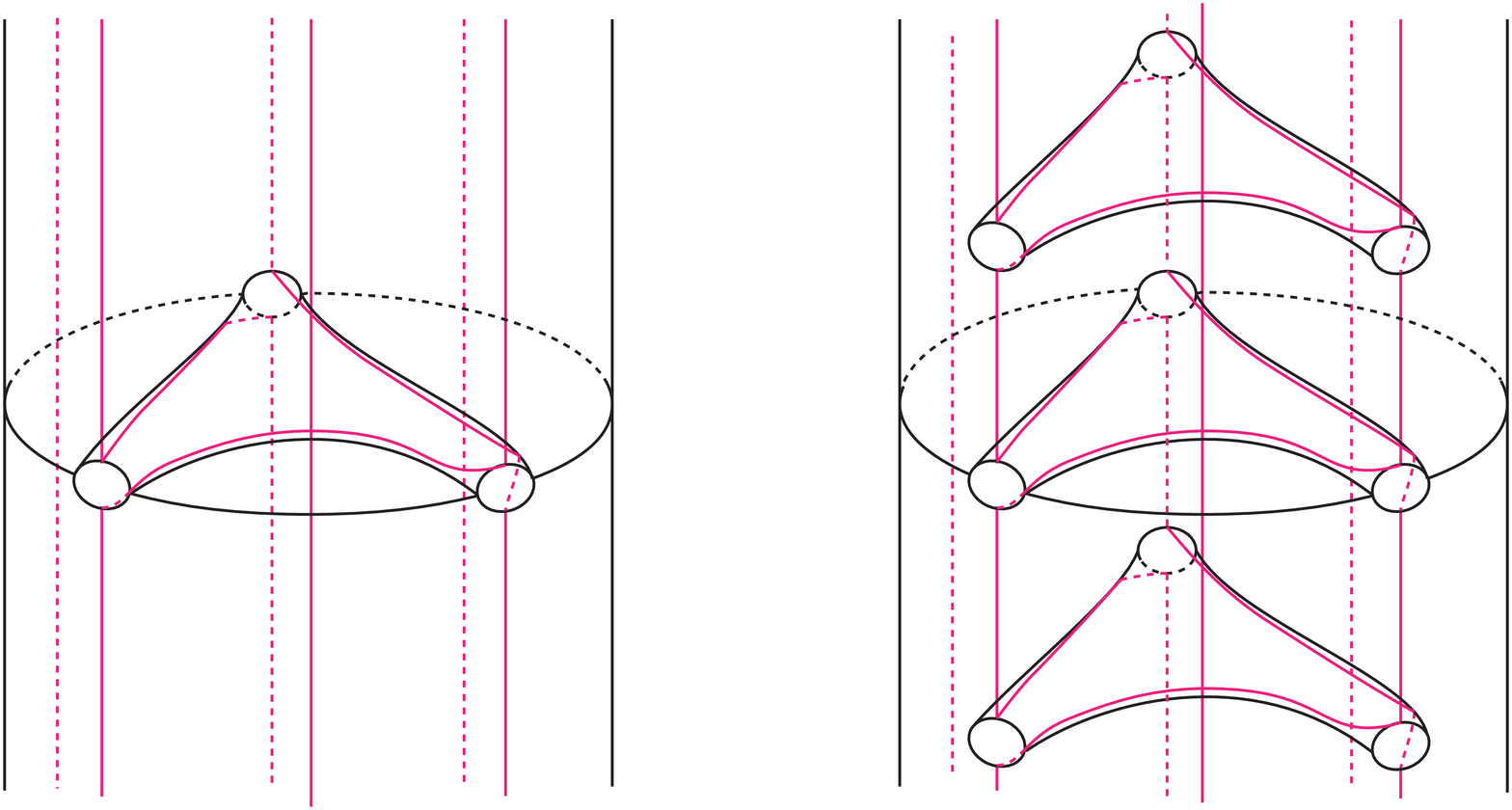}
\put(16,21.25){\tiny $D_i$}
\end{overpic}
\caption{The left-hand diagram shows $N(L)$, and the right-hand
diagram shows $N(L\cup L_{-\varepsilon}\cup L_\varepsilon)$.}
\label{mult}
\end{figure}

Next we extend $L$ to a Legendrian skeleton $K$ for
$(M,\Gamma,\xi)$, which satisfies the conditions of
Theorem~\ref{thm: decomposition}.  Before doing this, it would be
convenient to ``protect'' $T$ and $L$ as follows.  Since $T$ is
convex, there exists a $[-1,1]$-invariant neighborhood
$T\times[-1,1]$ of $T=T_0$. Let $T_{-\varepsilon}$, $T_0$ and
$T_\varepsilon$ be parallel copies of $T_0$, where
$T_t=T\times\{t\}$.  Then take parallel Legendrian graphs
$L_{i\varepsilon}$, $i=-1,0,1$, on $T_{i \varepsilon}$, and thicken
to obtain $N(L_{i \varepsilon})$. Here $L=L_0$. Now, apply the
technique of Theorem~\ref{thm: decomposition} to extend
$L_{-\varepsilon}\cup L_0\cup L_\varepsilon$ to $K$. We may assume
that the extension does not intersect the region between
$T_{-\varepsilon}$ and $T_{\varepsilon}$ and is disjoint from $
L_{-\varepsilon}\cup L_0\cup L_\varepsilon$. (This is because we can
use the compressing disks $D_i$, together with their translates
corresponding to $T_{\pm \varepsilon}$.)

Let $(S,R_+(\Gamma),h)$ be the partial open book decomposition of
$(M,\Gamma,\xi)$ corresponding to the decomposition into $N(K)$ and
$M-N(K)$.  The three connected components of
$S-\overline{R_+(\Gamma)}$ corresponding to $L_{i\varepsilon}$ will
be called $Q_{i\varepsilon}$. (Of course there may be other
connected components.) If $a$ is a properly embedded arc on $Q=Q_0$
(with endpoints away from $\Gamma$), let $a_{i\varepsilon}$ be the
corresponding copy on $Q_{i\varepsilon}$. Then $h(a)$ is a
concatenation $a'a_\varepsilon a''$, where $a'$ and $a''$ are arcs
which switch levels $Q\leftrightarrow Q_\varepsilon$ as given in
Figure~\ref{regionP}. The proof that $h(a)$ is indeed as described
is similar to the computation of Example 3 in Section~\ref{section:
examples}.

\begin{figure}[ht]
\begin{overpic}[width=7cm]{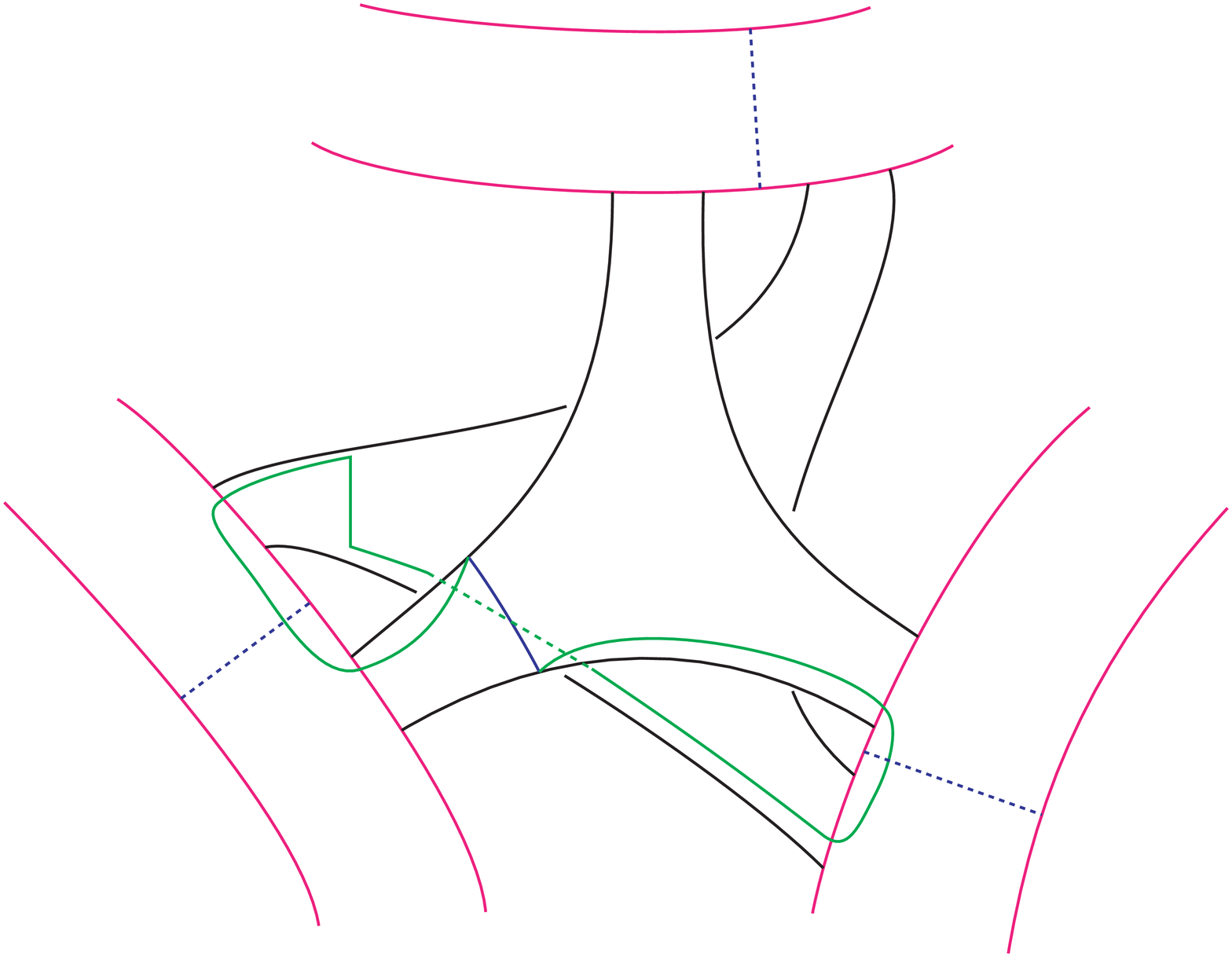}
\put(51,40){\tiny $Q$} \put(40,30){\tiny $a$} \put(63,50){\tiny
$Q_\varepsilon$} \put(29.5,35){\tiny $a_\varepsilon$}
\put(55,27){\tiny $a'$} \put(16,27){\tiny $a''$}
\end{overpic}
\caption{} \label{regionP}
\end{figure}

Next, a page of the partial open book decomposition
$(S',R_+(\Gamma'),h')$ of $(M',\Gamma',\xi')$ is given in
Figure~\ref{newopenbook}.  In particular, all of $Q$ becomes part of
the new $R_+(\Gamma')$ and we cut $R_+(\Gamma)\subset S$ along the
arcs $d_i$ that correspond to the $D_i$.  (If the $d_i$ are the
dotted blue arcs in Figure~\ref{regionP}, then $D_i=
d_i\times[-1,1]$.)

\begin{figure}[ht]
\begin{overpic}[width=7cm]{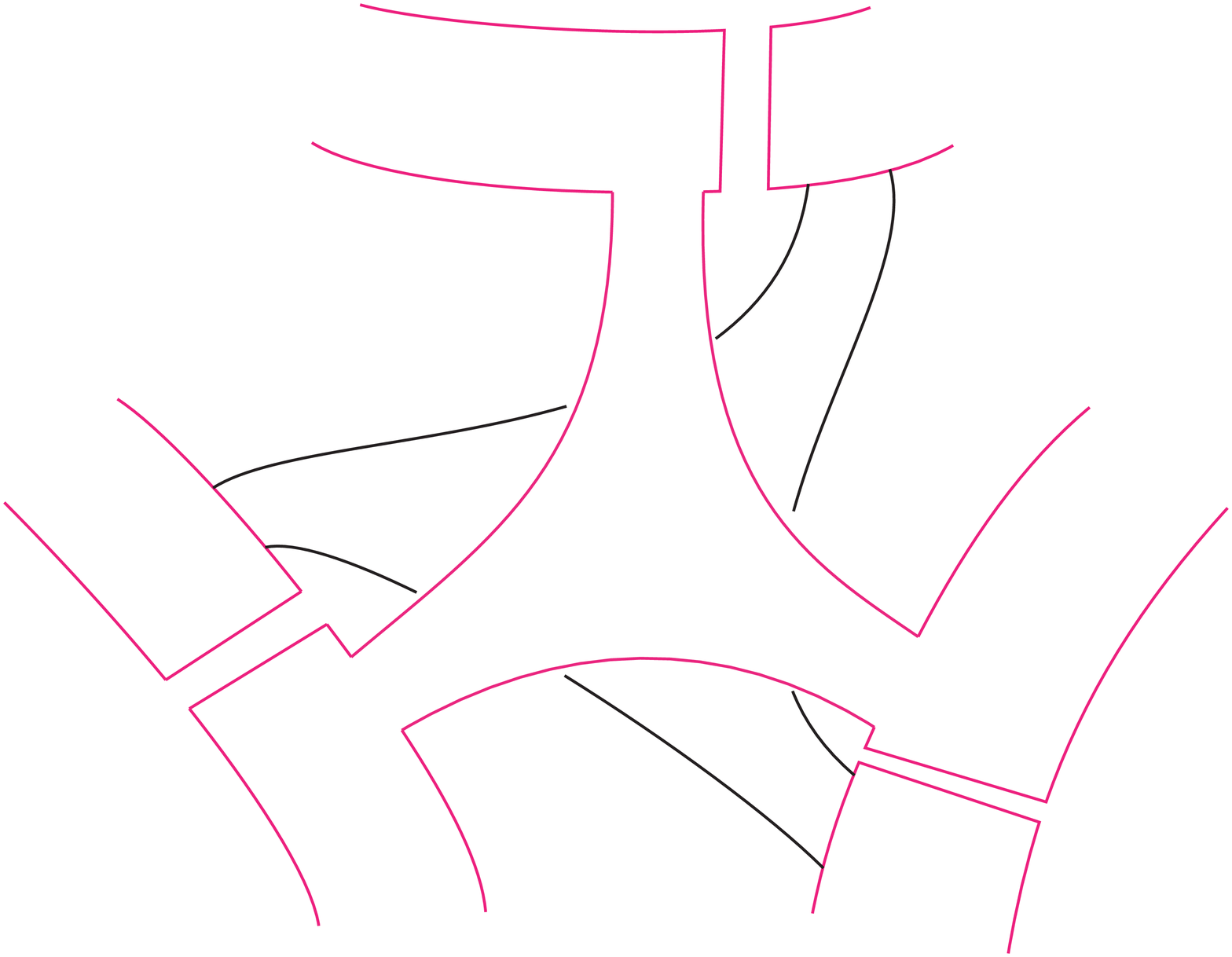}

\end{overpic}
\caption{} \label{newopenbook}
\end{figure}

\s\n {\bf Case 1: $Q$ has no genus.} Suppose $Q$ has no genus and
hence is a $2k$-gon, where $k$ of the edges are subarcs
$\gamma_0,\dots,\gamma_{k-1}$ of $\Gamma$ and $k$ of the edges are
subarcs of $\bdry S$. Assume the $\gamma_i$ are in counterclockwise
order about $\bdry Q$. Consider a basis for $Q$ consisting of $k-1$
arcs $a_1,\dots,a_{k-1}$, where $a_i$ is parallel to $\gamma_i$. Let
$b_i$ be the pushoffs of $a_i$ satisfying (1), (2) and (3) of
Section~\ref{section: defn} and let $x_i$ be the intersection of
$a_i$ and $b_i$ on $Q\times \{1\}$. Also let $\alpha_i=\bdry
(a_i\times[-1,1])$ and the $\beta_i=(b_i\times\{1\})\cup
(h(b_i)\times\{-1\})$, viewed on $S\times[-1,1]/_\sim$. (This
$[-1,1]$-coordinate is different from the one used previously for
the invariant neighborhood of $T$.) Complete the $\alpha_i$ and
$\beta_i$ into a compatible Heegaard decomposition for $(M,\Gamma)$
by adding $\alpha_i'$ and $\beta_i'$. Hence
$\alpha=(\alpha_1,\dots,\alpha_{k-1},\alpha_1',\dots,\alpha'_l)$ and
$\beta=(\beta_1,\dots,\beta_{k-1},\beta_1',\dots,\beta'_l)$.

\begin{claim}
The only $(k+l-1)$-tuples $\mathbf{y}$ of $\mathbb{T}_\alpha\cap
\mathbb{T}_\beta$ whose corresponding first Chern class
$c_1(\mathfrak{s}_\mathbf{y})$ evaluates maximally on $T$ have the
form $(y_1,\dots,y_{k-1},y'_1,\dots,y'_l)$, where $y_i$ is an
intersection of some $\alpha_{j_1}$ and $\beta_{j_2}$, and $y'_i$ is
an intersection of some $\alpha'_{j_1}$ and $\beta'_{j_2}$.
\end{claim}

\begin{proof}
We describe $T$ as it sits inside the partial open book
$(S,R_+(\Gamma),h)$. First isotop the $d_i$ above to $d_i'$ as given
in Figure~\ref{T}. Let $T'$ be the union of the $d_i'\times[-1,1]$
and the region $Q'\times\{1\}\subset Q\times\{1\}$ bounded by the
$d_i'\times\{1\}$ and the $\gamma_i$. Then isotop $T'$ to obtain $T$
so that $T\cap \Sigma$ is the union of $\mbox{(the shaded
region)}\times \{-1\}$ and $Q\times\{1\}$.  Here $\Sigma$ is the
Heegaard surface for $(S,R_+(\Gamma),h)$.
\begin{figure}[ht]
\begin{overpic}[width=7cm]{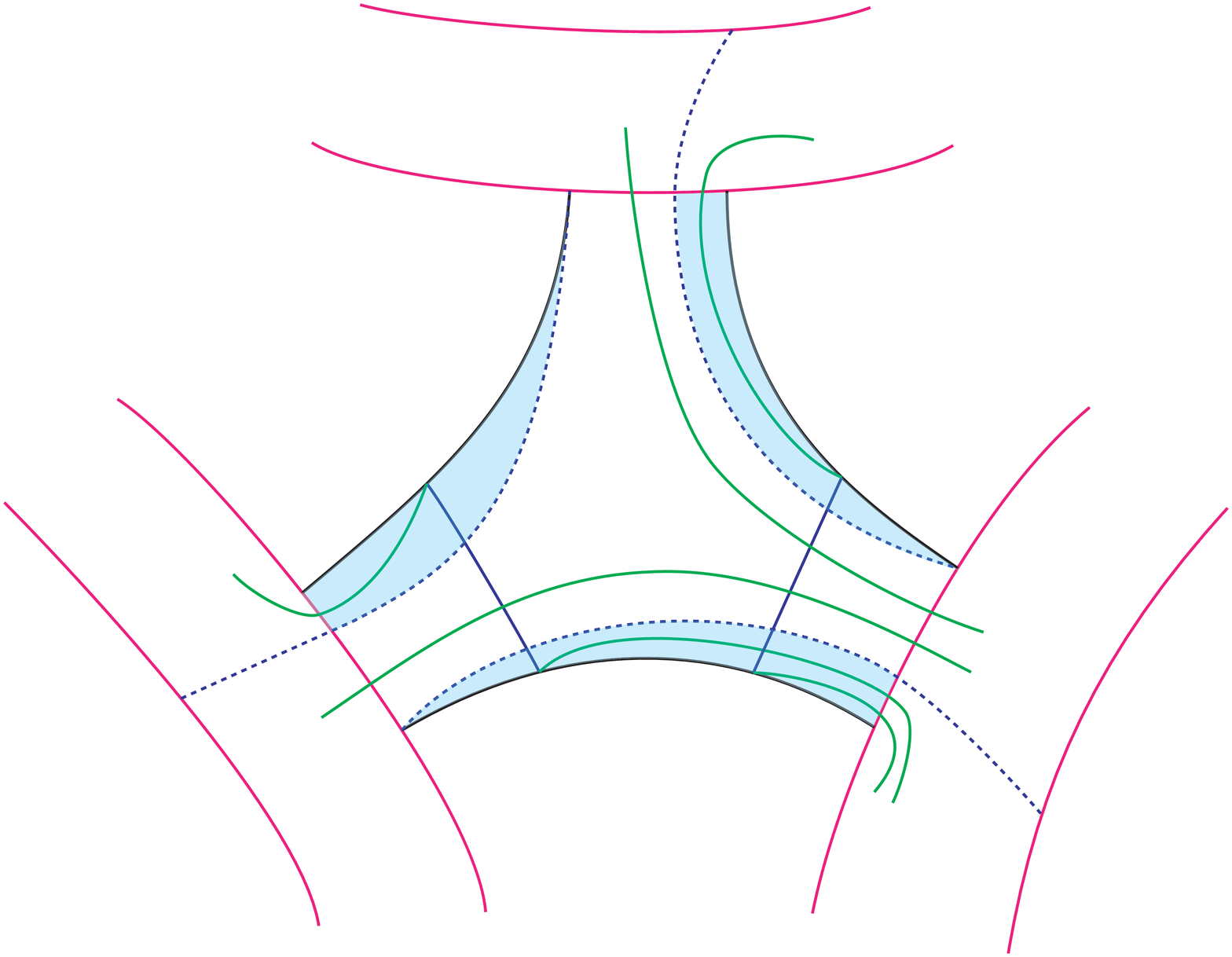}
\put(51,40){\tiny $Q$} \put(18,20){\tiny $d_i'$} \end{overpic}
\caption{} \label{T}
\end{figure}
Observe that the gradient flow line corresponding to $x_i\in
\alpha_i\cap \beta_i$ intersects $T$.  The same also holds for any
$y_i\in \alpha_{j_1}\cap \beta_{j_2}$, since $h(b_{j_2})\cap Q$ is
contained in the shaded region. On the other hand, the gradient flow
line for any intersection point of $\alpha_{j_1}$ and $\beta'_{j_2}$
does not intersect $T$, as $\beta'_{j_2}$ does not enter
$Q\times\{1\}$ or the shaded region due to the ``protection''
afforded by the $Q_{i\varepsilon}$. (The only $\beta'_{j_2}$ that
enter $Q\times\{-1\}$ come from $Q_{-\varepsilon}$.  In that case we
see in Figure~\ref{T} that they are represented by the green arcs
outside the shaded region.) Therefore, in order to maximize $\langle
c_1(\mathfrak{s}_\mathbf{y}), T\rangle$, the intersection point on
$\alpha_i$ must lie on $\beta_j$.  This forces $\alpha_i$ to be
paired with $\beta_j$ and $\alpha_i'$ to be paired with $\beta_j'$.
\end{proof}

Still assuming that $Q$ has no genus, consider $\mathbf{y}$ as in
the above claim.  The only $\beta_i$ which intersects $\alpha_1$ is
$\beta_1$, and their sole intersection is $x_1$. Hence $x_1$ occurs
in $\mathbf{y}$ --- this uses up $\alpha_1$ and $\beta_1$. Next, the
only $\beta_i$ besides $\beta_1$ which intersects $\alpha_2$ is
$\beta_2$. Continuing in this manner, $\{x_1,\dots, x_{k-1}\}\subset
\mathbf{y}$.  Hence, the inclusion
$$\Phi:\CF(\{\beta_1',\dots,\beta_l'\},\{\alpha_1',\dots,\alpha_l'\})\rightarrow
\CF(\beta,\alpha),$$
$$(y_1',\dots,y_l')\mapsto (x_1,\dots,x_{k-1},y_1',\dots,y_l'),$$
is as a direct summand of $\CF(\beta,\alpha)$.  Moreover, when
counting holomorphic disks from
$\mathbf{y'}=(x_1,\dots,x_{k-1},y_1',\dots,y_l')$ to
$\mathbf{y''}=(x_1,\dots,x_{k-1},y_1'',\dots,y_l'')$, the
positioning of $\Gamma$ and the $x_i$ implies that no subarcs of
$\alpha_i$ and $\beta_i$ can be used as the boundary of a nontrivial
holomorphic map.  Hence holomorphic disks from $\mathbf{y'}$ to
$\mathbf{y''}$ are in 1-1 correspondence with holomorphic disks from
$(y_1',\dots,y_l')$ to $(y_1'',\dots,y_l'')$.  This implies that
$\Phi$ is a chain map.  The homology of
$\CF(\{\beta_1',\dots,\beta_l'\},$ $\{\alpha_1',\dots,\alpha_l'\})$
is $\SFH(-M',-\Gamma')$, and image under the injection $\Phi$ is the
part of $\SFH(-M,-\Gamma)$ which is outer.  It is also easy to see
that $\EH(M',\Gamma',\xi')$ is mapped to $\EH(M,\Gamma,\xi)$ under
$\Phi$.

\begin{figure}[ht]
\begin{overpic}[width=10cm]{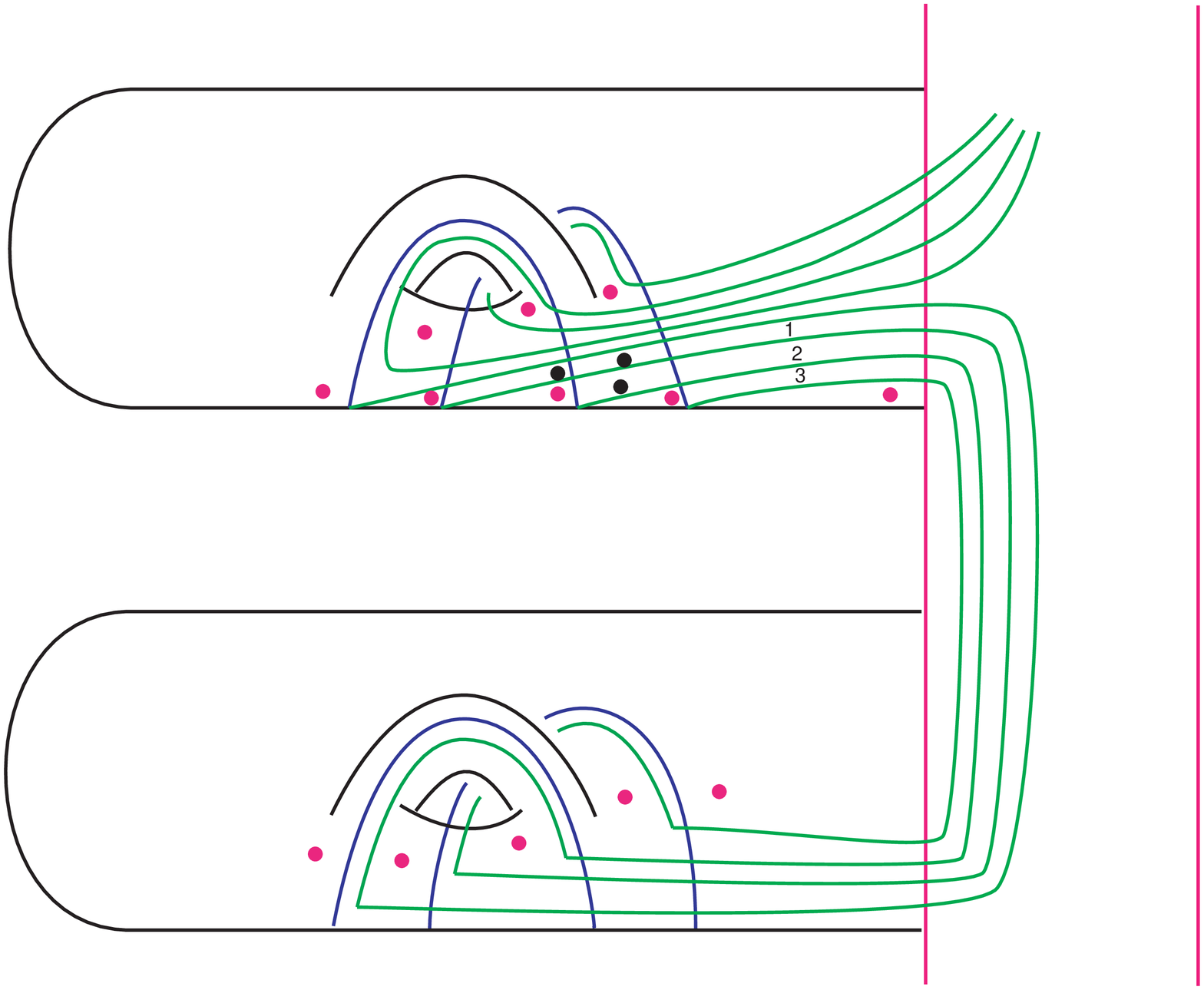}
\put(15,60){\tiny $Q$} \put(15,15){\tiny $Q_\varepsilon$}
\put(88,70){\tiny $\beta_j'$} \put(50,25){\tiny $\alpha_j'$}
\end{overpic} \caption{} \label{genus}
\end{figure}

\s\n {\bf Case 2: $Q$ has genus.}  In general, $Q$ has genus $g$ and
there are $k$ attaching arcs along $\Gamma$. For a general picture,
we can think of Figure~\ref{regionP} with added handles. The basis
for $Q$ can be chosen to consist of the same arcs
$a_1,\dots,a_{k-1}$ parallel to $\gamma_1,\dots,\gamma_{k-1}$ as in
Figure~\ref{regionP}, plus additional arcs $a_{k+2j-2}, a_{k+2j-1}$,
with one pair for each handle $j=1,\dots, g$. These additional arcs
can be taken to end on the subarc $\tau$ of $\bdry S$ between
$\gamma_{k-1}$ and $\gamma_0$, where $\gamma_{k-1}$,
$\tau$,$\gamma_0$ are in counterclockwise order about $\bdry Q$.
Figure~\ref{genus} depicts the case when $k=1, g=1$.

Let $a_1,\dots,a_m$, $m=k+2g-1$, be a basis for $Q$, and $b_i$,
$x_i$, $\alpha_i$, $\beta_i$, $i=1,\dots,m$,  and
$\alpha_j',\beta_j'$, $j=1,\dots,l$, be as before. The curves that
nontrivially intersect $Q\times\{-1\}$ are $\alpha_i$, $\beta_i$,
and curves $\beta_j'$ of the form $((b_i)_{-\varepsilon}
\times\{1\})\cup (h((b_i)_{-\varepsilon})\times\{-1\})$. Here
$(b_i)_{-\varepsilon}$ is an arc on $Q_{-\varepsilon}$ which is
parallel to $b_i\subset Q$. Moreover, the curves $\alpha_j'$ and
$\beta_j'$ which pass through $Q_{\varepsilon}\times\{-1\}$ can be
taken to be $\bdry((a_i)_\varepsilon\times [-1,1])$ and
$((b_i)_\varepsilon\times \{1\})\cup
(h((b_i)_\varepsilon)\times\{-1\})$, respectively.  This is due to
the protective layer between $T_{-\varepsilon}$ and
$T_{\varepsilon}$. See Figure~\ref{genus}. We can represent $T$ as
before, as the union of $\mbox{(the shaded region)}\times \{-1\}$
and $Q\times\{1\}$, and prove the analogous Claim by using this
description of $T$ to show that the $\alpha_{j_1}$ and $\beta_{j_2}$
must be paired to maximize the first Chern class on $T$. Hence we
only consider $\mathbf{y}=(y_1,\dots,y_m,y'_1,\dots,y'_l)$, where
$y_i$ is an intersection of some $\alpha_{j_1}$ and $\beta_{j_2}$,
and $y'_i$ is an intersection of some $\alpha'_{j_1}$ and
$\beta'_{j_2}$. As before, if $k>1$, then we argue that the only
$\beta_i$ which intersects $\alpha_1$ is $\beta_1$, the only
$\beta_i$ besides $\beta_1$ which intersects $\alpha_2$ is
$\beta_2$, etc., and that $\{x_1,\dots, x_{k-1}\}\subset
\mathbf{y}$.

Since the subset $\{x_1,\dots,x_{k-1}\}$ has no effect on what
follows, we assume that $k=1$.  Moreover, for simplicity, we assume
that $g=1$. (The higher genus case is not much more difficult.)
Consider the pair of curves $\alpha_1, \alpha_2$ (corresponding to
the handle) and the corresponding $\beta$ curves. Let $x_1, v_1$ be
the intersection points of $\alpha_1 \cap \beta_1$,  $v_2$ be the
intersection point of $\alpha_1 \cap \beta_2$, $w_1,w_3$ be the
intersection points of $\alpha_2 \cap \beta_1$, and $x_2,w_2$ be the
intersection points of $\alpha_2 \cap \beta_2$. See
Figure~\ref{genus2}.  Hence the summand $W\subset \CF(\beta,\alpha)$
corresponding to the generators which evaluate maximally on $T$ is
generated by $(y_1,y_2,\mathbf{y'})$, where $(y_1,y_2)$ is one of
$$(x_1,x_2), (x_1, w_2), (v_1,x_2), (v_1,w_2), (v_2,w_1), (v_2,w_3),$$
\begin{figure}[ht]
\s
\begin{overpic}[width=8cm]{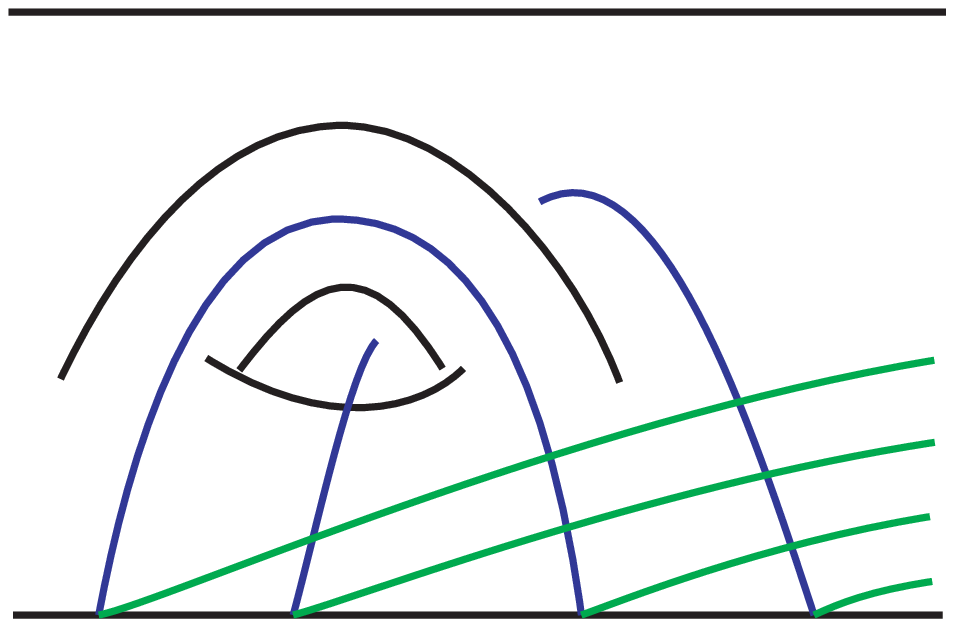}
\put(8.3,-4){\tiny $x_1$} \put(28.8,-4){\tiny $x_2$}
\put(59,-4){\tiny $x_1$} \put(84,-4){\tiny $x_2$}
\put(28.6,10){\tiny $u_1$} \put(53,18){\tiny $v_1$}
\put(55,11){\tiny $v_2$} \put(78,25){\tiny $w_1$}
\put(81,17.7){\tiny $w_2$} \put(83.5,10){\tiny $w_3$}
\end{overpic} \caption{} \label{genus2}
\end{figure}
and $\mathbf{y'}$ is of the type $(y_1',\dots,y_l')$.

Now we examine the holomorphic disks in $Sym^{m+l}(\Sigma)$ that
avoid $\bdry\Sigma\times Sym^{m+l-1}(\Sigma)$, in order to compute
the boundary maps. The realization of the holomorphic disk as a
holomorphic map $\overline{u}$ to $\Sigma$ cannot have ``mixed''
boundary components, i.e., each boundary component must consist
solely of subarcs of $\alpha_i$ and $\beta_j$, or consist solely of
subarcs of $\alpha'_i$ and $\beta'_j$, since all generators of $W$
consist of intersections of $\alpha_j$ and $\beta_k$, and
intersections of $\alpha_j'$ and $\beta'_k$ only. Consider the
regions of $\Sigma-\cup_i \alpha_i-\cup_i \beta_i-\cup_j\alpha_j'
-\cup_j\beta_j'$ which nontrivially intersect $\Gamma$ --- they are
indicated by red dots in Figure~\ref{genus}.  As a consequence, all
the regions in $Q\times\{-1\}$ besides the regions with black dots
and the regions $1$, $2$, and $3$, have multiplicity zero. Now,
regions $1$, $2$, and $3$ have $\alpha'_j$ curves on the boundary.
Hence, if the multiplicity is nonzero, the image of $\overline{u}$
must extend across those curves, in the end engulfing the red-dotted
regions in $Q_\varepsilon\times \{-1\}$.  Therefore, the
multiplicities of regions $1$, $2$, and $3$ are forced to be zero.

The following are the only possible holomorphic disks: \be
\item quadrilaterals in $Q\times\{-1\}$ (the regions with black
dots in Figure~\ref{genus});
\item holomorphic disks from $\mathbf{y'}$ to
$\mathbf{y''}$ that do not enter $Q\times\{-1\}$. \ee The
quadrilaterals give rise to boundary maps:
$$\bdry (u_1,v_2)=(x_2,v_1),$$
$$\bdry (v_1,w_2)=(v_2,w_1),$$
$$\bdry (v_2,w_3)=(x_1,w_2),$$
$$\bdry (x_1,x_2)=\bdry (x_2,v_1)=\bdry (v_2,w_1)=\bdry
(x_1,w_2)=0,$$ where we are referring to the Heegaard diagram with
the $\alpha'_i$, $\beta'_i$ erased. If $\mathbf{y_i}$ denotes a
linear combination of tuples of type $(y_1',\dots,y_l')$ (by slight
abuse of notation), then
$$\bdry ((x_1,x_2,\mathbf{y_1})+(u_1,v_2,\mathbf{y_2})+(x_2,v_1,\mathbf{y_3})+\dots)=0$$
implies
$$(x_1,x_2,\bdry \mathbf{y_1}) +(x_2,v_1,\mathbf{y_2})+(u_1,v_2,\bdry \mathbf{y_2})+(x_2,v_1,\bdry
\mathbf{y_3})+\dots=0.$$  Here $\bdry \mathbf{y_i}$ refers to the
boundary map for the Heegaard diagram with
$\alpha_1,\alpha_2,\beta_1,\beta_2$ erased. This implies that $\bdry
\mathbf{y_1}=\bdry \mathbf{y_2}=0$, and
$\bdry\mathbf{y_3}=\mathbf{y_2}$.  Hence,
$$(x_1,x_2,\mathbf{y_1})+(u_1,v_2,\mathbf{y_2})+(x_2,v_1,\mathbf{y_3})+\dots=
(x_1,x_2,\mathbf{y_1})+\bdry(u_1,v_2,\mathbf{y_3})+\dots,$$ and the
inclusion
$$\Phi:\CF(\{\beta_1',\dots,\beta_l'\},\{\alpha_1',\dots,\alpha_l'\})\rightarrow
\CF(\beta,\alpha),$$
$$(y_1',\dots,y_l')\mapsto (x_1,x_2,y_1',\dots,y_l'),$$
induces an isomorphism of $\SFH(-M',-\Gamma')$ onto a direct summand
of $\SFH(-M,-\Gamma)$, as before.
\end{proof}

\s\n {\em Acknowledgements.}  We thank Andr\'as Juh\'asz for
patiently explaining sutured Floer homology and Andr\'as Stipsicz
for helpful discussions.  We also thank John Etnyre for the
discussions which led to the reformulation of the contact invariant
in \cite{HKM2}.

\end{document}